\input amstex
\documentstyle{amsppt}
\magnification=\magstep1

\pageheight{9.0truein}
\pagewidth{6.5truein}
\TagsOnLeft
\NoBlackBoxes
\loadbold

\def\leftmat{\left[ \matrix}
\def\rightmat{\endmatrix \right]}

\def\AA{{\Bbb A}}
\def\CC{{\Bbb C}}

\def\QQ{{\Bbb Q}}
\def\RR{{\Bbb R}}
\def\ZZ{{\Bbb Z}}

\def\L{{\Cal L}}
\def\O{{\Cal O}}
\def\P{{\Cal P}}
\def\gg{{\frak g}}
\def\mm{{\frak m}}

\def\ahat{\widehat{a}}
\def\bhat{\widehat{b}}
\def\abar{\overline{a}}
\def\bbar{\overline{b}}

\def\hbar{\overline{h}}
\def\kbar{\overline{k}}

\def\Lbar{\overline{L}}
\def\Mbar{\overline{M}}
\def\Nbar{\overline{N}}

\def\Atil{\widetilde{A}}

\def\mun{\underline{m}}

\def\Oq{\O_q}
\def\kx{k^\times}
\def\bfp{\boldsymbol{p}}
\def\bfq{\boldsymbol{q}}
\def\chr{\operatorname{char}}
\def\gr{\operatorname{gr}}
\def\Sg{S(\gg)}
\def\Ug{U(\gg)}
\def\Ogs{\O(\gg^*)}
\def\Cinfty{C^\infty}
\def\gl{\frak{gl}}
\def\ad{\operatorname{ad}}
\def\Dx{\operatorname{Dx}}
\def\Dxbar{\operatorname{\overline{Dx}}}
\def\Dxtil{\operatorname{\widetilde{Dx}}}
\def\maxspec{\operatorname{maxspec}}
\def\Gspec{\operatorname{spec^{G}}}
\def\gspec{\operatorname{spec^{\gg}}}
\def\Obfq{\O_{\bfq}}
\def\CCx{\CC^\times}

\def\Pspec{\operatorname{P{.}spec}}
\def\Pprim{\operatorname{P{.}prim}}
\def\spec{\operatorname{spec}}
\def\prim{\operatorname{prim}}
\def\Fract{\operatorname{Fract}}

\def\AvMV {{\bf 1}}
\def\BGR{{\bf 2}}
\def\Bro{{\bf 3}}
\def\BrGd{{\bf 4}}
\def\BGY{{\bf 5}}
\def\BrGr{{\bf 6}}
\def\CG{{\bf 7}}
\def\CoMc{{\bf 8}}
\def\Dix{{\bf 9}}
\def\Dbook{{\bf 10}}
\def\FL{{\bf 11}}
\def\GMur{{\bf 12}}
\def\GAth{{\bf 13}}
\def\GLau{{\bf 14}}
\def\qaff{{\bf 15}}
\def\GLDM{{\bf 16}}
\def\GLquo{{\bf 17}}
\def\GLet{{\bf 18}}
\def\GoYa{{\bf 19}}
\def\Hay{{\bf 20}}
\def\Hod{{\bf 21}}
\def\HoLeA{{\bf 22}}
\def\HoLeB{{\bf 23}}
\def\HLT{{\bf 24}}
\def\Hrt{{\bf 25}}
\def\Ing{{\bf 26}}
\def\Jos{{\bf 27}}
\def\Kam{{\bf 28}}
\def\Kirnilp{{\bf 29}}
\def\Kirsurv{{\bf 30}}
\def\KirC{{\bf 31}}
\def\KlSc{{\bf 32}}
\def\KoSo{{\bf 33}}
\def\Mar{{\bf 34}}
\def\Mat{{\bf 35}}
\def\McPe{{\bf 36}}
\def\Nou{{\bf 37}}
\def\Ohone{{\bf 38}}
\def\Ohtwo{{\bf 39}}
\def\Ohthree{{\bf 40}}
\def\OPS{{\bf 41}}
\def\Pol{{\bf 42}}
\def\RTF{{\bf 43}}
\def\Shaf{{\bf 44}}
\def\SoiA{{\bf 45}}
\def\SoiB{{\bf 46}}
\def\Str{{\bf 47}}
\def\Tak{{\bf 48}}
\def\Taktwo{{\bf 49}}
\def\TY{{\bf 50}}
\def\VS{{\bf 51}}
\def\Vnc{{\bf 52}}
\def\Van{{\bf 53}}
\def\Wein{{\bf 54}}

\topmatter

\title Semiclassical limits of quantized coordinate rings\endtitle

\author K. R. Goodearl\endauthor

\dedicatory Dedicated to S. K. Jain on the occasion of his 70th birthday
\enddedicatory

\address Department of Mathematics, University of California,
Santa Barbara, CA 93106, USA\endaddress
\email goodearl\@math.ucsb.edu\endemail

\subjclassyear{2000}
\subjclass 16W35; 16D60, 17B63, 20G42
\endsubjclass

\keywords Quantized coordinate ring, semiclassical limit, Poisson
algebra, symplectic leaf, symplectic core, Dixmier map
\endkeywords

\abstract This paper offers an expository account of some ideas, methods,
and conjectures concerning quantized coordinate rings and their
semiclassical limits, with a particular focus on primitive ideal spaces.
The semiclassical limit of a family of quantized coordinate rings of an
affine algebraic variety $V$ consists of the classical coordinate ring
$\O(V)$ equipped with an associated Poisson structure. Conjectured
relationships between primitive ideals of a generic quantized coordinate
ring $A$ and symplectic leaves in $V$ (relative to a semiclassical
limit Poisson structure on $\O(V)$) are discussed, as are breakdowns in
the connections when the symplectic leaves are not algebraic. This
prompts replacement of the differential-geometric concept of symplectic
leaves with the algebraic concept of symplectic cores, and a reformulated
conjecture is proposed: The primitive spectrum of $A$ should be
homeomorphic to the space of symplectic cores in $V$, and
to the Poisson-primitive spectrum of $\O(V)$. Various examples, including
both quantized coordinate rings and enveloping algebras of solvable Lie
algebras, are analyzed to support the choice of symplectic cores to
replace symplectic leaves.
\endabstract

\endtopmatter

\document

\head 0. Introduction \endhead

By now, the ``Cheshire cat'' description of quantum groups is
well known -- a quantum group is not a group at all, but something that
remains when a group has faded away, leaving an algebra of functions
behind. The appropriate functions depend on which category of group is
under investigation. We concentrate here on (affine) algebraic groups $G$,
on which the natural functions of interest are the polynomial functions.
These constitute the classical coordinate ring of $G$, which we denote
$\O(G)$. (The group structure on $G$ induces a Hopf algebra structure on
$\O(G)$, but we shall not make use of that.) A {\it quantized
coordinate ring\/} of $G$ is, informally, a deformation of $\O(G)$, in
the sense that it is an algebra with a set of generators patterned after
those in $\O(G)$, but with a new multiplication that is typically
noncommutative. Examples and references will be given in Section 1.
We do not address the question of what properties are required to qualify
an algebra as a quantized coordinate ring -- this remains a fundamental
open problem. Quantized coordinate rings have also been defined for a
number of algebraic varieties other than algebraic groups, and our
discussion will incorporate them as well.

Many parallels have been found between the structures
of quantized and classical coordinate rings, and general principles for
organizing and predicting such parallels are needed. The present paper
concentrates on a circle of ideas and results focussed on ideal
structure, particularly spaces of prime or primitive ideals. The
theme/principle we follow, based on much previous work, can be stated this
way:
\roster
\item"$\bullet$" {\sl The primitive ideals of a suitably generic
quantized coordinate ring of an algebraic variety $V$ should match
subsets of $V$ in some partition defined through the geometry of $V$ and
a Poisson structure obtained from a semiclassical limit process.}
\endroster
Many of  the terms just mentioned require explanations, which we will
give over the course of the paper. Here we just mention that, in the
above statement, ``generic'' refers to the assumption that suitable
parameters in the construction of the quantized coordinate ring should be
non-roots of unity. 
\medskip

To begin the story (omitting many definitions and details), we refer to
the results of Soibelman and Vaksman \cite{\VS, \SoiA, \SoiB},
who studied the ``standard'' generic quantized coordinate rings of
simple compact Lie groups $K$. They established a bijection between the
irreducible *-representations of $K$  (on Hilbert spaces) and the
symplectic leaves in
$K$ (relative to a Poisson structure arising from the
quantization). This amounts to a linkage between primitive ideals and
symplectic leaves, a relationship which is a key ingredient of the Orbit
Method from Lie theory. Informed by this principle, and inspired by the
work of Soibelman and Vaksman, Hodges and Levasseur conjectured that
similar bijections should exist for semisimple complex algebraic groups
\cite{\HoLeA}. The case of $SL_2(\CC)$ being easy \cite{\HoLeA, Appendix},
they first verified the conjecture for $SL_3(\CC)$ \cite{op\. cit.}, and
then for
$SL_n(\CC)$ \cite{\HoLeB}. In later work with Toro \cite{\HLT}, they
verified it for connected semisimple groups. In light of these
achievements, it is natural to pose this conjecture for other classes of
generic quantized coordinate rings. (It is easily seen that the above
conjecture cannot hold for non-generic quantized coordinate rings. In
such cases, the quantized coordinate rings are usually finitely generated
modules over their centers, and they have far more primitive ideals than
can be matched to symplectic leaves.)

In the specific cases just mentioned, the symplectic leaves turn out to
be algebraic, in the sense that they are locally closed in the Zariski
topology. Hodges, Levasseur, and Toro pointed out in \cite{\HLT} that
symplectic leaves need not be algebraic for Poisson structures arising
from multiparameter quantizations, and that the above conjecture cannot
be expected to hold in such cases. We argue that it should not be
surprising that the concept of symplectic leaves, which comes from
differential geometry, is not always well suited for algebraic problems.
Thus, symplectic leaves should be replaced by more algebraically defined
objects. The notion of symplectic cores introduced by Brown and Gordon
\cite{\BrGr} fills the role well, up to the present state of
knowledge; we will give evidence to buttress this statement.
\medskip

Our aim here is to present an account of the above story, with
introductions to and discussions of the relevant concepts. In particular,
the tour will pass through way stations such as {\it quantized coordinate
rings\/}, {\it semiclassical limits\/}, {\it Poisson structures\/}, {\it
symplectic leaves\/}, the {\it Orbit Method\/}, {\it symplectic
cores\/}, and the {\it Dixmier map\/}. By the end of the tour, we will be
in purely algebraic territory, where we can formulate a conjecture that
does not require any differential geometry (i.e., symplectic leaves).
Namely:
\roster
\item"$\bullet$" {\sl If $A$ is a generic quantized coordinate ring
of an affine algebraic variety $V$ over an algebraically closed field of
characteristic zero, and if $V$ is given the Poisson structure arising
from an appropriate semiclassical limit, then the spaces of
primitive ideals in $A$ and symplectic cores in $V$, with their
respective Zariski topologies, are homeomorphic.}
\endroster
A parallel conjecture relates the prime and primitive spectra of $A$ to
the spaces of Poisson prime and Poisson-primitive ideals in $\O(V)$.
\medskip

Fix a base field $k$ throughout the paper; all algebras mentioned will be
unital $k$-algebras. This field can be general at first, but then we will
require it to have characteristic zero, and/or be algebraically closed.
When discussing symplectic leaves, we restrict $k$ to $\RR$ or $\CC$.

\head 1. Quantized coordinate rings \endhead

We begin by recalling two basic examples, to clarify the idea that a
quantized coordinate ring of an algebraic group (or variety) is, loosely
speaking, a deformation of the classical coordinate ring. References to
many other examples are given in
\S\S1.2, 1.4, 1.5.

\definition{1.1\. Quantum $SL_2$} Recall that the group $SL_2(k)$ is a
closed subvariety of the variety of $2\times2$ matrices over $k$, defined
by the single equation ``determinant = $1$''. The coordinate ring of the
matrix variety is naturally realized as a polynomial ring in four
variables $X_{ij}$, corresponding to the functions that pick out the four
entries of the matrices. The coordinate ring of $SL_2(k)$ can thus be
described as follows:
$$\O(SL_2(k)) = k[X_{11},X_{12},X_{21},X_{22}]/ \langle X_{11}X_{22} -
X_{12}X_{21} - 1 \rangle.$$
To ``quantize'' this coordinate ring, we replace the commutative
multiplication by a noncommutative one, parametrized by a nonzero scalar
$q$, as below. The reasons for this particular choice of relations
will not be given here; see \cite{\BrGd, \S\S I.1.6, I.1.8}, for instance,
for a discussion. 

Given a choice of scalar $q\in\kx$, the ``standard'' one-parameter {\it
quantized coordinate ring of $SL_2(k)$\/} is the $k$-algebra
$\Oq(SL_2(k))$ presented by generators $X_{11}$, $X_{12}$, $X_{21}$,
$X_{22}$ and the following relations:
$$\alignedat2 X_{11}X_{12} &= qX_{12}X_{11}  &X_{11}X_{21} &=
qX_{21}X_{11}\\
X_{12}X_{22} &= qX_{22}X_{12}  &X_{21}X_{22} &= qX_{22}X_{21}\\
X_{12}X_{21} &= X_{21}X_{12}  &\qquad\qquad X_{11}X_{22}-X_{22}X_{11} &=
(q-q^{-1}) X_{12}X_{21}\\
X_{11}X_{22} &- qX_{12}X_{21} = 1\,. \endalignedat$$
The case when $q=1$ is special: The first six relations then reduce to
saying that the generators $X_{ij}$ commute with each other, the last
reduces to the defining relation for the variety $SL_2(k)$, and so the
algebra
$\O_1(SL_2(k))$ is just the classical coordinate ring. We write this,
very informally, as
$$\O(SL_2(k)) = {\underset{q\rightarrow1} \to \lim}\ \Oq(SL_2(k));$$
it is our first example of a ``semiclassical limit''.
\enddefinition

\definition{1.2\. Quantum matrices, quantum $SL_n$ and $GL_n$} The
pattern indicated in \S1.1 extends to definitions of ``standard'' single
parameter quantized coordinate rings $\O_q(M_n(k))$, $\O_q(SL_n(k))$, and
$\O_q(GL_n(k))$ for all positive integers $n$. Multiparameter versions, which we label in the form
$\O_{\lambda,\bfp}(-)$, have also been defined. Generators and relations
for these algebras may be found, for instance, in \cite{\GMur,
\S\S1.2--1.4;
\BrGd, \S\S I.2.2--I.2.4}.
\enddefinition

\definition{1.3\. Quantum affine spaces} The coordinate ring of affine
$n$-space over $k$ is the polynomial algebra in $n$ indeterminates,
and the most basic quantization is obtained by replacing
commutativity $(xy=yx)$ with {\it $q$-commutativity\/}: $xy=qyx$. Thus,
the ``standard'' one-parameter {\it quantized coordinate ring of
$k^n$\/}, relative to a choice of scalar $q\in\kx$, is the $k$-algebra 
$$\Oq(k^n)= k\langle x_1,\dots,x_n \mid x_ix_j= qx_jx_i \text{\ for\ }
1\le i<j\le n \rangle.$$
The multiparameter version of this algebra requires an $n\times n$ matrix
of nonzero scalars, $\bfq= (q_{ij})$, which is {\it multiplicatively
antisymmetric\/} in the sense that $q_{ii}=1$ and $q_{ji}= q_{ij}^{-1}$
for all $i$, $j$. The {\it multiparameter quantized coordinate ring of
$k^n$\/} corresponding to a choice of $\bfq$ is the $k$-algebra
$$\Obfq(k^n)= k\langle x_1,\dots,x_n \mid x_ix_j= q_{ij}x_jx_i \text{\
for all\ } i,j \rangle.$$

In the one-parameter case, we can write
$\O(k^n)= {\underset{q\rightarrow1} \to \lim}\ \Oq(k^n)$
in the same sense as above. For the multiparameter case, we imagine a
limit in which all $q_{ij} \rightarrow 1$.
\enddefinition

\definition{1.4\. Quantized coordinate rings of semisimple groups} The
single parameter versions of these Hopf algebras, which we denote
$\Oq(G)$, were first defined for semisimple algebraic groups $G$ of
classical type (types A, B, C, D) via generators and relations, by
Faddeev, Reshetikhin, and Takhtadjan
\cite{\RTF} and Takeuchi \cite{\Tak}. A detailed development (done for
$k=\CC$, but the pattern is the same over other fields) can be found in
\cite{\KlSc, Chapter 9}. In most of the more recent literature, $\Oq(G)$ is
defined as a restricted Hopf dual of the quantized enveloping algebra of
the Lie algebra of $G$ (e.g., see
\cite{\BrGd, Chapter I.7}). This is a more uniform approach, which also
covers groups of exceptional type. That the two approaches yield the same
Hopf algebras in the classical cases was established by Hayashi
\cite{\Hay} and Takeuchi \cite{\Tak} (see \cite{\KlSc, Theorem 11.22}). 

The single parameter algebras $\Oq(G)$ constitute the ``standard''
quantized coordinate rings of semisimple groups. Multiparameter versions,
which we label $\O_{q,p}(G)$, were introduced by Hodges, Levasseur, and
Toro \cite{\HLT}.
\enddefinition

\definition{1.5\. Additional examples} Quantized coordinate rings, both
single- and multiparameter, have been defined for many algebraic
varieties, such as algebraic tori, toric varieties, and versions of
affine spaces related to classical groups of types B, C, D. For a general
survey, see \cite{\GMur, Section 1}.
Quantized toric varieties were introduced in \cite{\Ing} (see also
\cite{\GLquo, \GLDM}). A family of iterated skew polynomial algebras
covering multiparameter quantized euclidean and symplectic spaces was
introduced by Oh
\cite{\Ohone} and extended by Horton \cite{\Hrt} (see also \cite{\GLau,
\S2.5} for the odd-dimensional euclidean case). Among other algebras
that have been studied in the literature, we mention quantized
coordinate rings for varieties of antisymmetric matrices \cite{\Str} and
varieties of symmetric matrices \cite{\Nou, \Kam}.
\enddefinition

\definition{1.6\. Limits of families of algebras} The semiclassical
limits informally introduced in \S\S1.1, 1.3 are more properly viewed
in the framework of families of algebras. For example, the algebras
$\Oq(SL_2(k))$ are quotients of a single algebra over a Laurent
polynomial ring $k[t^{\pm1}]$, namely the algebra $A$ given by generators 
$X_{11}$, $X_{12}$, $X_{21}$,
$X_{22}$ and relations as in \S1.1, but with $q$ replaced by $t$:
$$\alignedat2 X_{11}X_{12} &= tX_{12}X_{11}  &X_{11}X_{21} &=
tX_{21}X_{11}\\
X_{12}X_{22} &= tX_{22}X_{12}  &X_{21}X_{22} &= tX_{22}X_{21}\\
X_{12}X_{21} &= X_{21}X_{12}  &\qquad\qquad X_{11}X_{22}-X_{22}X_{11} &=
(t-t^{-1}) X_{12}X_{21}\\
X_{11}X_{22} &- tX_{12}X_{21} = 1\,. \endalignedat  \tag 1.6a$$
For each $q\in\kx$, there is a natural identification $A/(t-q)A \equiv
\Oq(SL_2(k))$. The ``limit as $q\rightarrow 1$'' is then simply the case
$q=1$ of these identifications:  $A/(t-1)A \equiv
\O(SL_2(k))$.

Similarly, if we take
$$B= k[t^{\pm1}]\langle x_1,\dots,x_n \mid x_ix_j= tx_jx_i \text{\ for\ }
1\le i<j\le n \rangle,   \tag 1.6b$$
then $B/(t-q)B \equiv \Oq(k^n)$ for all $q\in\kx$, and
$${\underset{q\rightarrow1} \to \lim}\ \Oq(k^n) = B/(t-1)B \equiv
\O(k^n).$$
The multiparameter algebras $\Obfq(k^n)$ can, likewise, be set up as
common quotients of an algebra over a Laurent polynomial ring
$k[t_{ij}^{-1} \mid 1\le i<j\le n]$. However, for purposes such as
obtaining Poisson structures on semiclassical limits, we need to be able
to exhibit the $\Obfq(k^n)$ as quotients of $k[t^{\pm1}]$-algebras.
There are many ways to do this; we will discuss some in \S2.3.
\enddefinition

\definition{1.7\. An older example: the Weyl algebra} Weyl defined the
algebra we now call the {\it first Weyl algebra\/} as
$$\CC\langle x,y \mid xy-yx= \hslash i \rangle,$$
where $\hslash$ is Planck's constant and $i= \sqrt{-1}$. Physicists often
use the term ``classical limit'' to denote the transition from a quantum
mechanical system to a classical one by letting Planck's constant go to
zero. The fact that ${\underset{\hslash\rightarrow0} \to \lim}$ of
the above algebra is the polynomial ring $\CC[x,y]$ is one
instance of this point of view.

To relate this semiclassical limit to the ones above, take $k=\CC$ and
take the scalar $q$ in quantized coordinate rings to be $e^{\hslash}$.
Then $\hslash\rightarrow0$ corresponds to $q\rightarrow1$. In many
constructions, particularly the C*-algebra quantum groups corresponding
to compact Lie groups, the parameter $q$ is either written directly in
the form $e^{\hslash}$ or is taken to be a nonnegative real number, with
calculations involving $e^{\hslash}$ used for motivation.
\enddefinition

\head 2. Semiclassical limit constructions \endhead

In the context of quantized coordinate rings, semiclassical limits are
constructed via quotients of algebras over Laurent polynomial rings, as in
\S1.6. A different version, using associated graded rings, is
needed in other arenas, particularly for enveloping algebras of Lie
algebras. We describe both constructions in this section.

\definition{2.1\. Semiclassical limits: commutative fibre version} Let
$k[h]$ be a polynomial algebra, with the indeterminate named $h$ as a
reminder of Planck's constant. Suppose that $A$ is a torsionfree
$k[h]$-algebra, and that $A/hA$ is commutative. Since $A$ is then a flat
$k[h]$-module, the family of factor algebras $\bigl( A/(h-\alpha)A
\bigr)_{\alpha\in k}$ (or, for short, $A$ itself) is called a {\it
flat family\/} of $k$-algebras, and $A/hA$ is viewed as the ``limit'' of
the family. It may happen that some of the algebras $A/(h-\alpha)A$
collapse to zero or are otherwise not desirable. If so, it is natural to
treat $A$ as an algebra over a localization of $k[h]$ (cf\. Example
2.2(c), for instance). We will usually not do this explicitly.

An immediate question is, what kind of information about the algebras
$A/(h-\alpha)A$ is contained in this limit? Observe that, because of the
commutativity of $A/hA$, all additive commutators $[a,b]= ab-ba$ in $A$
are divisible by $h$. Moreover, division by $h$ is unique, since $A$ is
torsionfree as a $k[h]$-module. Hence, we obtain a well defined binary
operation $\frac1h[-,-]$ on $A$. This operation enjoys four key
properties:
\roster
\item Bilinearity;
\item Antisymmetry;
\item The Jacobi identity (thus $A$, equipped with
$\frac1h[-,-]$, is a Lie algebra over $k$);
\item The {\it Leibniz identities\/}, that is, the product
rule (for derivatives) in each variable: $\frac1h[a,bc]= \bigl(
\frac1h[a,b] \bigr)c + b
\bigl( \frac1h[a,c] \bigr)$ for all $a,b,c\in A$, and similarly for
$\frac1h[bc,a]$.
\endroster
Operations satisfying properties (1)--(4) are called {\it Poisson
brackets\/}.

The above Poisson bracket on $A$ induces, uniquely, a Poisson bracket on
$A/hA$, which we denote $\{-,-\}$. Thus, writing overbars to denote
cosets modulo $hA$, we have
$$\{\abar,\bbar\}= \overline{\tfrac1h[a,b]}$$
for $a,b\in A$. The commutative algebra $A/hA$, equipped with this
Poisson bracket, is called the {\it semiclassical limit\/} of the family
$\bigl( A/(h-\alpha)A
\bigr)_{\alpha\in k}$. Loosely speaking, the Poisson bracket on the
semiclassical limit records a ``first-order impression'' of the
commutators in $A$ and in the algebras $A/(h-\alpha)A$.
\enddefinition

\definition{2.2\. Examples} {\bf (a)} Fit the one-parameter quantum affine
spaces $\Oq(k^n)$ into the $k[t^{\pm1}]$-algebra $B$ of (1.6b), and
set $h=t-1$. Then $B$ represents a flat family of $k[h]$-algebras, with
$B/hB$ commutative. We identify $B/hB$ with the polynomial ring
$k[x_1,\dots,x_n]$ and compute the resulting Poisson bracket on the
indeterminates as follows.  For
$1\le i<j\le n$, we have
$[x_i,x_j]= hx_jx_i$ in
$B$, and hence
$$\{x_i,x_j\}= x_ix_j$$
in $k[x_1,\dots,x_n]$. Because of the Leibniz
identities, the above information determines this Poisson bracket
uniquely. It may be described in full as follows:
$$\{f,g\}= \sum_{1\le i<j\le n} x_ix_j \left( \frac{\partial f}{\partial
x_i} \frac{\partial g}{\partial x_j} - \frac{\partial g}{\partial
x_i} \frac{\partial f}{\partial x_j} \right)$$
for all $f,g\in k[x_1,\dots,x_n]$.

{\bf (b)} Take $A= k[h]\langle x,y\mid xy-yx= h \rangle$. Then
$A/(h-\alpha)A \cong A_1(k)$ for all nonzero $\alpha\in k$, while $A/hA$
can be identified with the polynomial ring $k[x,y]$. In this case,
the semiclassical limit Poisson bracket on $k[x,y]$ satisfies (and is
determined by)
$$\{x,y\} = 1.$$

{\bf (c)} The family $\bigl( \Oq(SL_2(k)) \bigr)_{q\in\kx}$ fits into the 
$k[t^{\pm1}]$-algebra $A$ with generators $X_{11}$, $X_{12}$, $X_{21}$,
$X_{22}$ and relations (1.6a). This is a flat family over
$k[h]$, where $h=t-1$. Since $t$ is invertible in $A$, the specialization
$A/(h+1)A$ is zero, corresponding to the fact that $A$ is actually a
torsionfree (even free) algebra over the localization $k[h][(h+1)^{-1}]$.
Here the semiclassical limit is the classical coordinate ring
$\O(SL_2(k))$, equipped with the Poisson bracket satisfying
$$\xalignat2
\{X_{11},X_{12}\} &= X_{11}X_{12}  &\{X_{11},X_{21}\} &= X_{11}X_{21}\\
\{X_{12},X_{22}\} &= X_{12}X_{22}  &\{X_{21},X_{22}\} &= X_{21}X_{22}\\
\{X_{12},X_{21}\} &= 0  &\{X_{11},X_{22}\} &= 2X_{12}X_{21} \,.
\endxalignat$$
\enddefinition

\definition{2.3\. Multiparameter examples} To obtain a semiclassical
limit -- with Poisson bracket -- for a multiparameter family of algebras,
we convert to a single parameter family and apply the construction of
\S2.1. The procedure is clear when the parameters involved are integer
powers of a single parameter. For example, consider the algebras
$\Obfq(k^n)$ where $\bfq= (q^{a_{ij}})$ for $q\in\kx$ and an
antisymmetric integer matrix $(a_{ij})$. Then define
$$A= k[t^{\pm1}] \langle x_1,\dots,x_n \mid x_ix_j= t^{a_{ij}}x_jx_i
\text{\ for all\ } i,j\rangle,$$
which is a torsionfree $k[t-1]$-algebra with $A/(t-1)A$ commutative. The
semiclassical limit is the polynomial algebra $k[x_1,\dots,x_n]$,
equipped with the Poisson bracket satisfying
$$\{x_i,x_j\}= a_{ij}x_ix_j$$
for all $i$, $j$.

More general parameters can be dealt with by various means. A simple but
ad hoc method to handle any $\Obfq(k^n)$ is via the algebra
$$A= k[h] \langle x_1,\dots,x_n \mid x_ix_j= \bigl( 1+(q_{ij}-1)h \bigr)  
x_jx_i \text{\ for\ } 1\le i<j\le n\rangle,$$
which is set up so that $A/(h-1)A \cong \Obfq(k^n)$ and $A/hA \cong
k[x_1,\dots,x_n]$. This yields a Poisson bracket satisfying $\{x_i,x_j\}=
(q_{ij}-1)x_ix_j$ for all $i,j$.

A variant of the previous procedure, involving quadratic rather than
linear polynomials in $h$, is used in \cite{\GLet} to construct
semiclassical limits for which the conjecture sketched in the Introduction
 applies to the generic multiparameter quantum affine spaces
$\Obfq(k^n)$.

Inverse to the construction of semiclassical limits is the problem of
{\it quantization\/}: trying to represent a given algebra supporting a
Poisson bracket as a semiclassical limit of a suitable family of
algebras. We will not discuss this problem except to indicate a solution
for the case of homogeneous quadratic Poisson brackets on polynomial
rings. Namely, suppose we have a polynomial algebra $k[x_1,\dots,x_n]$,
equipped with a Poisson bracket such that $\{x_i,x_j\}=
\alpha_{ij}x_ix_j$ for all $i$, $j$, where $(\alpha_{ij})$ is an
antisymmetric matrix of scalars over $k$. In place of ad hoc procedures
such as the one sketched above, it is natural, assuming that $\chr k=0$,
to use power series. In this case, set
$\text{exp}(\alpha_{ij})= \sum_{n=0}^\infty \frac1{n!} \alpha_{ij}^n h^n
\in k[[h]]$ for all $i$, $j$, and form the $k[[h]]$-algebra
$$A= k[[h]] \langle x_1,\dots,x_n \mid x_ix_j=
\text{exp}(\alpha_{ij})x_jx_i \text{\ for all\ } i,j\rangle.$$
The semiclassical limit algebra is $k[x_1,\dots,x_n]$, and its
Poisson bracket is the original one. 

Analogous $k[[h]]$-algebra constructions are given for commonly studied
families of quantized coordinate rings of skew polynomial type in
\cite{\GLau, Section 2}.
\enddefinition

\definition{2.4\. Semiclassical limits: filtered/graded version} Suppose
that $A$ is a $\ZZ$-filtered $k$-algebra, say with filtration
$(A_n)_{n\in\ZZ}$. Thus, the $A_n$ are $k$-subspaces of $A$, with
$A_m\subseteq A_n$ when $m\le n$, such that $A_mA_n\subseteq A_{m+n}$ for
all $m$, $n$. We will assume that the filtration is {\it exhaustive\/},
that is, that $\bigcup_{n\in\ZZ} A_n = A$. Note that we must have $1\in
A_0$; thus, $A_0$ is a unital subalgebra of $A$. Finally, let $\gr A=
\bigoplus_{n\in\ZZ}
\gr_n A$ be the associated graded algebra, where $\gr_n A= A_n/A_{n-1}$.

Now assume that
$\gr A$ is commutative. Homogeneous elements $a\in \gr_m A$ and
$b\in \gr_n A$ can be lifted to elements $\ahat\in A_m$ and $\bhat\in
A_n$, and since $\gr A$ is commutative, the commutator $[\ahat,\bhat]$
must lie in $A_{m+n-1}$. We then set $\{a,b\}$ equal to the coset of
$[\ahat,\bhat]$ in $\gr_{m+n-1} A$. It is an easy exercise, left to the
reader, to verify that $\{a,b\}$ is well defined, and that the extension
of $\{-,-\}$ to sums of homogeneous elements defines a Poisson bracket on
$\gr A$. The commutative algebra $\gr A$, equipped with this Poisson
bracket, is called the {\it semiclassical limit\/} of $A$.

More generally, assume there is an integer $d<0$ such that
$[A_m,\, A_n] \subseteq A_{m+n+d}$
for all $m,n\in\ZZ$. This assumption of course forces $\gr A$ to be
commutative. Modify the definition above by setting
$\{a,b\}$ equal to the coset of $[\ahat,\bhat]$ in $\gr_{m+n+d} A$, for
$a\in \gr_m A$ and $b\in \gr_n A$. This recipe again produces a well
defined Poisson bracket on $\gr A$ \cite{\Mar, Lemma 2.7}.
\enddefinition

\definition{2.5\. Bridging the two constructions} The semiclassical limit
of a $\ZZ$-filtered algebra $A$ constructed in \S2.4 can also be obtained
by applying the construction of \S2.1 to an auxiliary algebra, namely the
{\it Rees ring\/} 
$$\Atil := \sum_{n\in\ZZ} A_nh^n \subseteq A[h^{\pm1}],$$
where $A[h^{\pm1}]$ is a Laurent polynomial ring over $A$. Since $1\in
A_0$, the polynomial algebra $k[h]$ is a subalgebra of $A$, and we note
that $\Atil$ is a torsionfree $k[h]$-algebra. (It is not a
$k[h^{\pm1}]$-algebra unless $A_{-1}= A_0$, in which case all $A_n=A_0$.)
On one hand,
$\Atil/(h-1)\Atil \cong A$. On the other, $\Atil/h\Atil \cong \gr A$,
because $h\Atil= \sum_{n\in\ZZ} A_{n-1}h^n$. Thus, if $\gr A$ is
commutative, we have a Poisson bracket $\frac1h[-,-]$ on $\Atil$, which
induces a Poisson bracket $\{-,-\}_1$ on $\gr A$ as in \S2.1. This bracket
concides with the Poisson bracket $\{-,-\}_4$ constructed in \S2.4, as
follows.

Start with $a\in \gr_m A$ and
$b\in \gr_n A$, and lift these elements to $\ahat\in A_m$ and $\bhat\in
A_n$. With respect to the natural epimorphism $\pi: \Atil \rightarrow \gr
A$, the elements $a$ and $b$ lift to $\ahat h^m,\bhat h^n \in \Atil$.
Hence,  
$$\{a,b\}_1= \pi \bigl( \frac1h [\ahat
h^m, \bhat h^n] \bigr)= \pi( [\ahat,\bhat] h^{m+n-1})= [\ahat,\bhat]+
A_{m+n-2} = \{a,b\}_4 \,.$$
Therefore $\{-,-\}_1= \{-,-\}_4$.
\enddefinition

\definition{2.6\. Example: enveloping algebras} Let $\gg$ be a finite
dimensional Lie algebra over $k$, and put the standard (nonnegative)
filtration on the enveloping algebra $\Ug$, so that $\Ug_0=k$ and $\Ug_1=
k+\gg$, while $\Ug_n= \Ug_1^n$ for $n>1$. The associated graded
algebra is commutative, and is naturally identified with the symmetric
algebra $\Sg$ of the vector space $\gg$. In particular, we use the same
symbol to denote an element of $\gg$ and its coset in $\gr_1\Ug=
\Sg_1$. Then $\Sg$ is the semiclassical limit of $\Ug$, equipped with the
Poisson bracket satisfying
$$\{e,f\}= [e,f]$$
for all $e,f\in \gg$, where $[e,f]$ denotes the Lie product in
$\gg$. The above formula determines $\{-,-\}$ uniquely, since $\gg$
generates $\Sg$.

Now view the dual space $\gg^*$ as an algebraic variety, namely the
affine space $\AA^{\dim\gg}$. The coordinate ring $\Ogs$ is a
polynomial algebra over $k$ in $\dim\gg$ indeterminates, as is $\Sg$.
There is a canonical isomorphism
$$\theta: \Sg @>{\ \cong\ }>> \Ogs  \tag 2.6$$
which sends each $e\in\gg$ to the polynomial function on
$\gg^*$ given by evaluation at $e$, that is, $\theta(e)(\alpha)=
\alpha(e)$ for $\alpha\in \gg^*$. (This isomorphism is often treated as an
identification of the algebras $\Sg$ and $\Ogs$.) Via $\theta$,
the Poisson bracket on $\Sg$ obtained from the semiclassical limit
process above carries over to a Poisson bracket on $\Ogs$, known as the
{\it Kirillov-Kostant-Souriau Poisson bracket\/}.

If $\{e_1,\dots,e_n\}$ is a basis for
$\gg$, then $\Sg= k[e_1,\dots,e_n]$ and $\theta$ sends the $e_i$ to
indeterminates $x_i$ such that
$\Ogs= k[x_1,\dots,x_n]$. An explicit description of the KKS Poisson
bracket on $\Ogs$ can be obtained in terms of the structure constants of
$\gg$, as follows. These constants are scalars $c^l_{ij} \in k$ such
that $[e_i,e_j]= \sum_l c^l_{ij}e_l$ for all $i$, $j$. Since $\{e_i,e_j\}=
[e_i,e_j]$ in $\Sg$, an application of $\theta$ yields $\{x_i,x_j\}=
\sum_l c^l_{ij}x_l$ for all $i$, $j$. It follows that
$$\{p,q\}= \sum_{i,j,l} c^l_{ij}x_l \frac{\partial p}{\partial x_i}
\frac{\partial q}{\partial x_j}$$
for $p,q\in \Ogs$ \cite{\CG,
Proposition 1.3.18}. To see this, just check that the displayed formula
determines a Poisson bracket on $\Ogs$ which agrees with the KKS bracket
on pairs of indeterminates.

The KKS Poisson bracket on $\Ogs$ can also be obtained by applying the
method of \S2.1 to the homogenization of $\Ug$, that is, the
$k[h]$-algebra $A$ with generating vector space $\gg$ and relations
$ef-fe= h[e,f]$ for $e,f\in\gg$ (where $[e,f]$ again denotes the Lie
product in $\gg$). Here $A/hA\cong \Sg \cong \Ogs$ and
$A/(h-\lambda)A\cong \Ug$ for all $\lambda\in\kx$.
\enddefinition

\head 3. Symplectic leaves \endhead

We introduce symplectic leaves first in the context of Poisson
 manifolds, following the original definition of Weinstein \cite{\Wein},
and then we carry the concept over to complex affine Poisson varieties,
following Brown and Gordon \cite{\BrGr}. 

\definition{3.1\. Poisson algebras} We reiterate the general definition
from \S2.1: a {\it Poisson bracket\/} on a $k$-algebra $R$ is any
antisymmetric bilinear map
$R\times R\rightarrow R$ which satisfies the Jacobi and Leibniz
identities. Unless a special notation imposes itself, we denote all
Poisson brackets by curly braces: $\{-,-\}$.

A {\it Poisson algebra\/} over $k$ is
just a $k$-algebra $R$ equipped with a particular Poisson bracket. We
restrict our attention to commutative Poisson algebras in the present paper. As for the noncommutative case, Farkas and Letzter have shown that Poisson brackets essentially reduce to commutators \cite{\FL, Theorem 1.2}: If $R$ is a prime ring which is not
commutative, any Poisson bracket on $R$ is a multiple of the commutator
bracket by an element of the extended centroid of $R$.
\enddefinition

\definition{3.2\. Symplectic leaves in Poisson manifolds} Let $M$ be a
smooth manifold, and let $\Cinfty(M)$ denote the algebra of smooth
real-valued functions on $M$. (Some authors replace $\Cinfty(M)$ by the
algebra of smooth or analytic complex-valued functions.) A {\it Poisson
structure\/} on $M$ is a choice of Poisson bracket on $\Cinfty(M)$, so
that
$\Cinfty(M)$ becomes a Poisson algebra. A smooth manifold, together with
a choice of Poisson structure, is called a {\it Poisson manifold\/}.

Now assume that $M$ is a Poisson manifold. For each $f\in \Cinfty(M)$,
the map $X_f= \{f,-\}$ is a derivation on $\Cinfty(M)$ and thus a vector
field on $M$. Such vector fields are called {\it Hamiltonian vector
fields\/} (for the given Poisson structure), and the flows (or integral
curves) of Hamiltonian vector fields are known as {\it Hamiltonian
paths\/}. More specifically, a smooth path $\gamma: [0,1] \rightarrow M$
is Hamiltonian provided there is some $f\in \Cinfty(M)$ such that, at each
point $\gamma(t)$ along the path, the tangent vector $d\gamma/dt$ equals
$X_f|_{\gamma(t)}$. Since the change from a Hamiltonian
path following the flow of a vector field $X_f$ to one following a
different vector field $X_g$ need not be smooth, one must work with {\it
piecewise Hamiltonian paths\/}, i.e., finite concatenations of
Hamiltonian paths.

These paths determine an equivalence relation on $M$, points $p$ and $p'$
being equivalent if and only if there is a piecewise Hamiltonian path in
$M$ running from $p$ to $p'$. The resulting equivalence classes are called
{\it symplectic leaves\/}, and the partition of $M$ as the disjoint
union of its symplectic leaves is known as the {\it symplectic
foliation\/} of $M$.   
\enddefinition

\definition{3.3\. Poisson bivector fields} For many purposes, it is more
useful to record a Poisson structure in the form of a bivector field
rather than a Poisson bracket. In particular, this allows the most direct
definition of Poisson structures on non-affine algebraic varieties.

Let $M$ be a Poisson manifold. For a point $p\in M$, let $\mm_p$ denote
the maximal ideal of $\Cinfty(M)$ consisting of those functions that
vanish at
$p$. Evaluation of Poisson brackets at $p$ induces an antisymmetric
bilinear form $\pi_p$ on the cotangent space $\mm_p/\mm^2_p$, where
$$\pi_p(f+\mm_p^2,\, g+\mm_p^2)= \{f,g\}(p)$$
 for $f,g\in \mm_p$. Now $\pi_p$ acts in each variable as a linear map in
the dual space of $\mm_p/\mm_p^2$, that is, as a tangent vector to $M$ at
$p$. Since $\pi_p$ is antisymmetric, it is thus a {\it tangent
bivector\/} at $p$, namely an element of $T_p(M)\wedge T_p(M)$. The map
$\pi: p\mapsto \pi_p$ is a smooth global section of $\Lambda^2 T_M$, that
is, a {\it tangent bivector field\/} on $M$. To recover the Poisson
bracket on $\Cinfty(M)$ from the bivector field $\pi$, observe that
$\{f,g\}(p)= \{f-f(p),\, g-g(p)\}(p)$ for $f,g\in \Cinfty(M)$ and $p\in
M$, which we rewrite in the form
$$\{f,g\}(p)= \pi_p(df(p),dg(p)),  \tag  3.3$$
where $df(p)= f-f(p)+\mm_p^2 \in \mm_p/\mm_p^2$
and similarly for $dg(p)$.

Conversely, via (3.3) any tangent bivector field $\pi$ on $M$ induces
an antisymmetric bilinear map $\{-,-\}$ on $\Cinfty(M)$ satisfying the
Leibniz conditions. This is a Poisson bracket exactly when the Jacobi
identity is satisfied, which is equivalent to the vanishing of the {\it
Schouten bracket\/} $[\pi,\pi]$ (which we will not define here; see
\cite{\AvMV, p\. 44; \Van, 2nd\. ed., Remark 2.2(3)}, for instance). A
{\it Poisson bivector field\/} on
$M$ is any tangent bivector field $\pi$ for which $[\pi,\pi] =0$. As
indicated in the sketch above, Poisson brackets on $\Cinfty(M)$
correspond bijectively to Poisson bivector fields on $M$.
\enddefinition

\definition{3.4\. Poisson varieties} For any complex algebraic variety
$V$, the definition of a {\it Poisson bivector field\/} on $V$ can be
copied from \S3.3 -- it is any tangent bivector field $\pi$ on $V$ for
which $[\pi,\pi] =0$. In the context of algebraic geometry, however, the
map $\pi: V \rightarrow \Lambda^2 T_V$ is required to be a regular
function. Now one defines a {\it Poisson variety\/} to be a complex
algebraic variety equipped with a particular Poisson bivector field.
Associated concepts are defined by requiring compatibility with these
bivector fields. For example, a {\it Poisson morphism\/} from a Poisson
variety $(V,\pi)$ to a Poisson variety $(W,\pi')$ is a regular map $\phi:
V \rightarrow W$ such that $(T\phi\wedge T\phi)\pi= \pi'\phi$. A {\it
Poisson subvariety\/} of $V$ is a subvariety $X$ such that the inclusion
map $X\rightarrow V$ is a Poisson morphism.

If $V$ is an affine Poisson variety, the formula (3.3) defines a
Poisson bracket on $\O(V)$. Conversely, any Poisson bracket on $\O(V)$
induces a Poisson bivector field on $V$ as in \S3.3. Thus, affine Poisson
varieties can equally well be defined as complex affine varieties whose
coordinate rings are Poisson algebras. This point of view can be extended
to arbitrary varieties by defining a Poisson variety to be a complex
algebraic variety whose sheaf of regular functions is a sheaf of Poisson
algebras.
\enddefinition

\definition{3.5\. Smooth Poisson varieties as manifolds} In order to
define symplectic leaves in Poisson varieties, manifold structures are
needed. The fundamental result is that any smooth (i.e., nonsingular)
complex variety $V$ has a unique structure as a complex analytic manifold
(e.g., \cite{\Shaf, Chapter II, \S2.3}). This allows one to view $V$ as a
smooth manifold. If
$V$ is a Poisson variety, its chosen Poisson bivector field $\pi$ is
necessarily smooth (because it is regular), and so $V$ together with
$\pi$ becomes a Poisson manifold. One can achieve this result with
Poisson brackets as well, by showing that any Poisson bracket on $\O(V)$
extends uniquely to a Poisson bracket on the algebra of smooth complex
functions on $V$; taking real parts then yields a Poisson bracket on
$\Cinfty(V)$.

Given a smooth Poisson variety $V$, we view $V$ as a smooth manifold as
above, and define the {\it symplectic leaves\/} of $V$ to be the
symplectic leaves of the manifold $V$, defined as in \S3.2. 
\enddefinition

\definition{3.6\. Symplectic leaves in singular Poisson varieties} Let
$V$ be an arbitrary complex variety, and define the sequence of closed
subvarieties
$$V_0=V \supset V_1\supset \cdots\supset V_m =\varnothing,$$
where each $V_{i+1}$ is the singular locus of $V_i$. To build this chain,
recall first that the singular locus of a nonempty variety is a proper
closed subvariety. Since $V$ is a noetherian topological space, the chain
must eventually reach the empty set.

If $V$ is a Poisson variety, then $V_1$ is a Poisson subvariety
\cite{\Pol, Corollary 2.4}. By induction, all the $V_i$ are Poisson
subvarieties of $V$. Consequently, $V$ is (canonically) the disjoint
union of smooth locally closed Poisson subvarieties $Z_i :=
V_{i-1}\setminus V_i$. Following \cite{\BrGr, \S3.5}, we define the {\it
symplectic leaves\/} of $V$ to be the symplectic leaves of the various
$Z_i$, defined as in \S3.5.
\enddefinition

\definition{3.7\. Example} There is a known recipe, described in
\cite{\HoLeA, Appendix A}, for determining the symplectic leaves
in a semisimple complex algebraic group
$G$, relative to the Poisson structure arising from the ``standard
quantization'' of $G$. For
illustration, we present the case $G= SL_2(\CC)$; details are
given in \cite{\HoLeA, Theorem B.2.1}. The Poisson bracket on $\O(G)$ is
described in \S2.2(c) above. The symplectic leaves in $G$ are as follows:
\roster
\item"$\bullet$" the singletons $\left\{ \leftmat \alpha&0\\ 0&\alpha^{-1}
\rightmat \right\}$, for $\alpha\in \CCx$; \smallskip
\item"$\bullet$" the sets $\left\{ \leftmat \alpha&0\\ \gamma&\alpha^{-1}
\rightmat \biggm| \alpha,\gamma \in \CCx \right\}$ and  $\left\{
\leftmat \alpha&\beta\\ 0&\alpha^{-1}
\rightmat \biggm| \alpha,\beta \in \CCx \right\}$; \smallskip
\item"$\bullet$" the sets  $\left\{ \leftmat \alpha&\beta\\ \gamma&\delta
\rightmat \in G \biggm| \beta=\lambda\gamma \ne 0 \right\}$, for $\lambda \in
\CCx$.
\endroster
\enddefinition

\definition{3.8\. Example} The standard example of a
non-algebraic solvable Lie algebra is a $3$-dim\-en\-sional complex Lie
algebra $\gg$ with basis
$\{e_1,e_2,e_3\}$ such that
$$\xalignat3 [e_1,e_2] &= e_2  &[e_1,e_3] &= \alpha e_3  &[e_2, e_3] &= 0
\endxalignat$$ 
for some $\alpha \in \RR\setminus \QQ$. Write $\Ogs=
\CC[x_1,x_2,x_3]$ following the notation of \S2.6. The KKS Poisson
structure on $\gg^*$ is given by the Poisson bracket on $\Ogs$ such that
$$\xalignat3 \{x_1,x_2\} &= x_2  &\{x_1,x_3\} &= \alpha x_3  &\{x_2,x_3\}
&= 0.
\endxalignat$$ 
As in \cite{\Van, 1st\. ed., Example II.2.37; 2nd\. ed., Example II.2.43},
the symplectic leaves in $\gg^*$ are the following sets:
\roster
\item"$\bullet$" the individual points on the $x_1$-axis;
\item"$\bullet$" the $x_1x_2$-plane minus the $x_1$-axis;
\item"$\bullet$" the $x_1x_3$-plane minus the $x_1$-axis;
\item"$\bullet$" the surfaces $(x_3= \lambda x_2^\alpha \ne 0)$ for
$\lambda\in\CCx$.
\endroster
Since $\alpha$ is irrational, the surfaces $(x_3= \lambda x_2^\alpha \ne
0)$ are not algebraic -- they are locally closed in the euclidean
topology but not in the Zariski topology.
\enddefinition

\head 4. The Orbit Method from Lie theory \endhead

\definition{4.1\. The Orbit Method} This term has been applied to a whole
complex of methods in the representation theory of Lie groups and Lie
algebras, and extended, as a guiding principle, to many other domains. To
quote Kirillov's survey article \cite{\Kirsurv}, 
\roster
\item"" {\sl The idea behind the orbit method is the unification of
harmonic analysis with symplectic geometry (and it can also be considered
as a part of the more general idea of the unification of mathematics and
physics). In fact, this is a post factum formulation. Historically, the
orbit method was proposed in \cite{\Kirnilp} for the description of the
unitary dual (i.e\. the set of equivalence classes of unitary irreducible
representations) of nilpotent Lie groups. It turned out that not only
this problem but all other principal questions of representation
theory---topological structure of the unitary dual, explicit description
of the restriction and induction functors, character formulae, etc.---can
be naturally answered in terms of coadjoint orbits.}
\endroster

In Lie theory, the relevant orbits are defined as follows. Recall that
if $G$ is a Lie group with Lie algebra $\gg$, then $G$ acts on $\gg$ by
the {\it adjoint action\/} and on $\gg^*$ by the {\it coadjoint
action\/}. The $G$-orbits of these actions are called {\it adjoint
orbits\/} and {\it coadjoint orbits\/}, respectively.
As a particular instance, the Orbit Method suggests that the primitive
ideals of the enveloping algebra of $\gg$, being
the kernels of the irreducible representations, should be related to the
coadjoint orbits in $\gg^*$. Kirillov's original work provided the best
such relationship -- a bijection -- when $\gg$ is nilpotent. There is
also a bijection in case $\gg$ is solvable, except that the coadjoint
orbits may have to be taken with respect to a different group than a Lie
group with Lie algebra $\gg$. We discuss this situation in Section 5.

To place the coadjoint orbits in a geometric setting, view $\gg^*$ as the
variety $\AA^{\dim\gg}$, as in \S2.6. We can then ask for a geometric
description of these orbits within $\gg^*$. The answer is a famous result
discovered independently by Kirillov, Kostant, and Souriau (see, e.g.,
\cite{\KirC, \S I.2.2, Theorem 2}):
\enddefinition

\proclaim{4.2\. Theorem} {\rm [Kirillov-Kostant-Souriau]} Let $G$ be a
Lie group and
$\gg$ its Lie algebra. Then the coadjoint orbits of
$G$ in $\gg^*$ are precisely the symplectic leaves for the KKS Poisson
structure. \endproclaim

\definition{4.3\. Example} Return to Example 3.8, and place the
$x_1x_2x_3$-coordinates of points of
$\gg^*$ in column vectors. We choose a Lie
group $G$ with Lie algebra $\gg$ as follows:
$$G := \left\{ \leftmat 1&u&v\\ 0&t&0\\ 0&0&t^\alpha
\rightmat \biggm| u,v\in\CC,\ t\in\CCx \right\}.$$
The coadjoint action of $G$ on $\gg^*$ can be identified with left
multiplication of matrices from $G$ on column vectors representing
points in $\gg^*$. One easily checks that the $G$-orbits are exactly the
symplectic leaves of $\gg^*$ identified in Example 3.8, as required by
Theorem 4.2.
\enddefinition

\definition{4.4\. A general principle} In situations outside Lie theory,
there may not be a suitable group action whose orbits play the role of
coadjoint orbits. Instead, taking account of Theorem 4.2, one focusses on
symplectic leaves. Restricting to the study of irreducible
representations and primitive ideals, one is led to a general principle
that we formulate as follows:
\roster
\item"" {\sl Given a noncommutative algebra $A$, relate the primitive
ideals of $A$ to the symplectic leaves corresponding to the Poisson
structure on some associated algebraic variety arising from a
semiclassical limit.}
\endroster
This loose phrasing is intended to give the flavor of ideas coming out of
the Orbit Method rather than to set up a precise recipe. Furthermore,
this principle already requires modification in the case of enveloping
algebras, and for general quantized coordinate rings. 

On the other hand, the principle is right on target for the generic
single parameter quantized coordinate rings $\Oq(G)$ of semisimple
complex algebraic groups $G$, as conjectured by Hodges and Levasseur in
\cite{\HoLeA, \S2.8, Conjecture 1}: there is a bijection between the set
of primitive ideals of $\Oq(G)$ and the set of symplectic leaves in $G$
(for the semiclassical limit Poisson structure). They verified this
conjecture for
$G= SL_2(\CC)$ and $G=SL_3(\CC)$ in \cite{\HoLeA, Corollary B.2.2,
Theorems 4.4.1, A.3.2}, and then for $G=SL_n(\CC)$ in \cite{\HoLeB,
Theorem 4.2 and following remarks}. For arbitrary connected, simply
connected, semisimple complex Lie groups $G$, Joseph proved in
\cite{\Jos, Theorem 9.2} that the primitive ideal space of
$\Oq(G)$ has the form conjectured by Hodges and Levasseur in the first
part of
\cite{\HoLeA, \S2.8, Conjecture 1}, but he did not address connections
with symplectic leaves. The full conjecture was established by Hodges,
Levasseur, and Toro \cite{\HLT, Theorems 1.8, 4.18, Corollary 4.5} for
connected semisimple complex Lie groups $G$. Their results also cover the
multiparameter algebra $\O_{q,p}(G)$ under suitable algebraicity
conditions on $p$.
\enddefinition

\definition{4.5\. Generic versus non-generic situations} As mentioned in
the introduction, the principle discussed in \S4.4 does not apply to
non-generic quantized coordinate rings, which typically have ``too many''
primitive ideals. The quantum plane provides the simplest illustration of
this difficulty, and of the differences between the generic and
non-generic cases. Take
$$A_q= \Oq(\CC^2)= \CC\langle x,y\mid xy=qyx\rangle,$$
where $q$ is an arbitrary nonzero scalar in $\CC$. By Example 2.2(a), the
semiclassical limit of the family $(A_q)_{q\in\CCx}$ is the
polynomial ring $k[x,y]$, equipped with the Poisson bracket such that
$\{x,y\}= xy$. It is easily checked that the corresponding symplectic
leaves in $\CC^2$ consist of
\roster
\item"$\bullet$" the individual points on the $x$- and $y$-axes;
\item"$\bullet$" the $xy$-plane minus  the $x$- and $y$-axes.
\endroster
If $q$ is not a root of unity, one similarly checks that the primitive
ideals of $A_q$ consist of
\roster
\item"$\bullet$" the maximal ideals $\langle x-\alpha,\,y\rangle$ and
$\langle x,\,y-\beta\rangle$, for $\alpha,\beta\in\CC$;
\item"$\bullet$" the zero ideal.
\endroster
(See \cite{\BrGd, Example II.7.2}, for instance, for details.) In this
case, there is a natural bijection between the set of primitive ideals of
$A_q$ and the set of symplectic leaves in $\CC^2$.

On the other hand, if $q$ is a primitive $l$-th root of unity, the center
of $A_q$ equals the polynomial ring $\CC[x^l,y^l]$, and $A_q$ is a
finitely generated $\CC[x^l,y^l]$-module. In this case, the primitive
ideals of $A_q$ are maximal ideals, and they are parametrized (up to
$l$-to-one) by the maximal ideals of $\CC[x^l,y^l]$. While the set of
primitive ideals of $A_q$ has the same cardinality as the set of
symplectic leaves in $\CC^2$, there is no natural bijection, and
certainly no homeomorphism if Zariski topologies are taken into account.

Such disparities occur in all the standard families of quantized
coordinate rings, and provide just one of many distinctions between the
generic and non-generic cases. We do not discuss the non-generic situation
further, and concentrate on generic algebras.
\enddefinition

\head 5. Limitations of the Orbit Method for solvable Lie algebras
\endhead

For a solvable finite dimensional
complex Lie algebra $\gg$, the primitive ideals of the enveloping algebra
$\Ug$ are parametrized by means of the famous {\it Dixmier map\/}. At
first glance, this is a successful instance of the Orbit Method, since the
Dixmier map induces a bijection from a set of orbits in $\gg^*$ onto the
set of primitive ideals of $\Ug$. However, the relevant orbits are not, in
general, those of the coadjoint action of a Lie group with Lie algebra
$\gg$. Instead, the following group is needed.

\definition{5.1\. The algebraic adjoint group} Let $\gg$ be a finite
dimensional complex Lie algebra. Treating
$\gg$ for a moment just as a vector space, we have the general linear
group
$GL(\gg)$ on $\gg$, which is a complex algebraic group whose Lie algebra
is the general linear Lie algebra $\gl(\gg)$. Any algebraic subgroup of
$GL(\gg)$ (i.e., any Zariski closed subgroup) has a Lie algebra which is
naturally contained in $\gl(\gg)$. The {\it algebraic adjoint group of
$\gg$\/} is the smallest algebraic subgroup $G\subseteq GL(\gg)$ whose
Lie algebra contains $\ad\gg= \{ \ad x\mid x\in \gg\}$ (cf\. \cite{\BGR,
\S12.2;
\TY, Definition 24.8.1}).

The natural action of $GL(\gg)$ on $\gg$ by linear automorphisms
restricts to an action of $G$ on $\gg$, the {\it adjoint action\/}. This,
in turn, induces a (left) action of $G$ on $\gg^*$, the {\it coadjoint
action\/}, under which
$$(g.\alpha)(x) = \alpha(g^{-1}.x)$$
for $g\in G$, $\alpha\in \gg^*$, and $x\in \gg$. The orbits of this
action, the {\it coadjoint orbits\/}, are collected in the set
$\gg^*/G$. We equip $\gg^*/G$ with the quotient topology induced from the
Zariski topology on
$\gg^*$, and thus refer to it as the {\it space\/} of coadjoint orbits.
\enddefinition

\definition{5.2\. Prime and primitive spectra} For any algebra $A$, we
denote the collection of all primitive ideals of $A$ by $\prim A$. This
set supports a Zariski topology, under which the closed sets are the sets
$V(I):= \{ P\in \prim A \mid P \supseteq I\}$ for ideals $I$ of $A$. We
treat
$\prim A$ as a topological space with this topology, and refer to it as
the {\it primitive spectrum\/} of $A$. The analogous process, applied to
the set of all prime ideals of $A$, results in the {\it prime spectrum\/}
of $A$, denoted $\spec A$. Since primitive ideals are prime, $\prim A
\subseteq
\spec A$. In fact, $\prim A$ is a subspace of $\spec A$, that is, its
topology coincides with the relative topology inherited from $\spec A$.
Finally, we shall need the subspace of $\spec A$ consisting of all the
maximal ideals of $A$. This is the {\it maximal ideal space\/} of $A$,
denoted $\maxspec A$.
\enddefinition

\definition{5.3\. The Dixmier map}  Let $\gg$ be a solvable finite
dimensional complex Lie algebra. Following \cite{\BGR, \S10.8}, we use
the name {\it Dixmier map\/} and the label $\Dx$ for the map
$$\Dx: \gg^*\longrightarrow \prim \Ug$$
introduced by Dixmier in \cite{\Dix}. We do not give the definition here,
but just refer to \cite{\BGR}. It turns out that this map is constant on
$G$-orbits, and so it induces a {\it factorized Dixmier
map\/}
$$\Dxbar: \gg^*/G \longrightarrow \prim \Ug$$
\cite{\BGR, \S12.4}. Work of Dixmier, Conze, Duflo, and Rentschler led to
the result that $\Dxbar$ is a continuous bijection \cite{\BGR, S\"atze
13.4, 15.1}. The conjecture that it is a homeomorphism was established
later by Mathieu \cite{\Mat, Theorem}, resulting in the following theorem:
\enddefinition

\proclaim{5.4\. Theorem} {\rm [Dixmier-Conze-Duflo-Rentschler-Mathieu]}
  Let $\gg$ be a solvable finite
dimensional complex Lie algebra, and $G$ its adjoint algebraic group.
Then  the factorized Dixmier map
$\Dxbar$ is a homeomorphism from $\gg^*/G$ onto $\prim\Ug$. \endproclaim

\definition{5.5\. Algebraic versus non-algebraic cases}
If $\gg$ is an {\it algebraic\/} Lie algebra, meaning that it is the Lie
algebra of some algebraic group, then the adjoint algebraic group $G$ is
a Lie group, and its coadjoint orbits in $\gg^*$ are the symplectic
leaves for the KKS Poisson structure, by Theorem 4.2. Otherwise, $G$ is
larger than the relevant Lie group, in the sense that its Lie algebra
properly contains $\ad\gg$. In this case, its coadjoint orbits are larger
too, typically larger than individual symplectic leaves. Our basic
example illustrates this behavior.
\enddefinition

\definition{5.6\. Example} Return to Example 3.8, and again place the
$x_1x_2x_3$-coordinates of points of
$\gg^*$ in column vectors. The adjoint algebraic group $G$, written
so as to act by left multiplication on column vectors, can be expressed as
$$G = \left\{ \leftmat 1&u&v\\ 0&t&0\\ 0&0&t' \rightmat
\biggm| u,v\in\CC,\ t,t'\in\CCx \right\}$$
\cite{\TY, \S24.8.4}. The coadjoint orbits of $G$ in $\gg^*$ are the
following sets:
\roster
\item"$\bullet$" the individual points on the $x_1$-axis;
\item"$\bullet$" the $x_1x_2$-plane minus the $x_1$-axis;
\item"$\bullet$" the $x_1x_3$-plane minus the $x_1$-axis;
\item"$\bullet$" $\gg^*$ minus the $x_1x_2$- and $x_1x_3$-planes.
\endroster
Comparing with Example 3.8, we see that the first three $G$-orbits are
symplectic leaves, while the fourth is not. However, the fourth is at
least a union of symplectic leaves. 

The fourth $G$-orbit above is Zariski dense in $\gg^*$, while the others
are not. Viewing these orbits as points in the orbit space $\gg^*/G$, we
find that $\gg^*/G$ has a unique dense point (i.e., a unique dense
singleton subset). By Theorem 5.4, the same holds for $\prim\Ug$.
(Translated into ideal theory, this means that there is one primitive
ideal of $\Ug$ which is contained in all other primitive ideals.) On the
other hand, all the surfaces
$(x_3= \lambda x_2^\alpha \ne 0)$ are Zariski dense in $\gg^*$, and so the
quotient topology on the space of symplectic leaves in $\gg^*$ has
uncountably many dense points. Therefore this space of symplectic leaves
cannot be homeomorphic to $\prim \Ug$.
\enddefinition

\head 6. Poisson ideal theory and symplectic cores \endhead

Since the concept of symplectic leaves is differential-geometric, it
should not be so surprising that it is not always suited to describe
answers to algebraic problems, as seen in the previous section.
Consequently, we look for an algebraic replacement. This is provided by
Brown and Gordon's notion of {\it symplectic cores\/}, which is described
via the ideal theory of Poisson algebras.

\definition{6.1\. Poisson prime ideals} Let $R$ be a (commutative) Poisson
algebra (recall \S3.1). 

A {\it Poisson ideal\/} of $R$ is any ideal $I$ of the ring $R$ which is
also a Lie ideal relative to $\{-,-\}$, that is, $\{R,I\} \subseteq I$.
Sums, products, and intersections of Poisson ideals are again
Poisson ideals. Whenever $I$ is a Poisson ideal of $R$, the Poisson
bracket on
$R$ induces a well defined Poisson bracket on $R/I$, so that $R/I$
becomes a Poisson algebra.

The {\it Poisson core\/} of an arbitrary ideal $J$ of $R$ is the largest
Poisson ideal contained in $J$. This exists and is unique, because it is
the sum of all Poisson ideals contained in $J$. We use $\P(J)$ to denote
the Poisson core of $J$.

A {\it Poisson-prime ideal\/} of $R$ is any proper Poisson ideal $P$ of
$R$ with the following property: whenever the product of Poisson ideals
$I$ and $J$ of $R$ is contained in $P$, one of $I$ or $J$ must be
contained in $P$. Obviously any prime Poisson ideal is Poisson-prime, but
the converse can fail in positive characteristic. As we shall see in a
moment, (Poisson-prime) is the same as (prime Poisson) when $R$ is
noetherian and $k$ has characteristic zero; in that case, we will drop the
hyphen and speak of {\it Poisson prime ideals\/}. Note also that if
$Q$ is an arbitrary prime ideal of
$R$, then
$\P(Q)$ is a Poisson-prime ideal.

The {\it Poisson-prime spectrum\/} of $R$, denoted $\Pspec R$, is the set
of all Poisson-prime ideals of $R$, equipped with the natural Zariski-type
topology, in which the closed sets are those of the form $V_P(I):= \{P\in
\Pspec R
\mid P \supseteq I\}$, for ideals $I$ of $R$. It suffices to consider
Poisson ideals in defining closed sets, since the ideal $I$ in the
definition of a closed set can be replaced by the Poisson ideal it
generates. (This observation is helpful in showing that finite
unions of closed sets are closed.) 
\enddefinition

\proclaim{6.2\. Lemma} Let $R$ be a Poisson $k$-algebra, where $\chr
k=0$. Then the Poisson core of every prime ideal of $R$ is prime, and all
minimal prime ideals of
$R$ are Poisson ideals. If
$R$ is noetherian, the Poisson-prime ideals of $R$ coincide with the prime
Poisson ideals.
\endproclaim

\demo{Proof} Commutativity is not needed for this result. The commutative case is covered, for instance, by
\cite{\GAth, Lemma 1.1}, and the general case is proved the same way. We
sketch the details for the reader's convenience.

The first conclusion is a consequence of \cite{\Dbook, Lemma 3.3.2}, and
the second follows.

Now assume that $R$ is noetherian, and let $P$ be a Poisson-prime ideal
of $R$. There exist prime ideals $Q_1,\dots,Q_t$ minimal over $P$ such
that $Q_1Q_2\cdots Q_t \subseteq P$. The minimal prime ideals $Q_i/P$ in
the Poisson algebra $R/P$ must be Poisson ideals by what has been proved
so far, and hence the $Q_i$ are Poisson ideals of $R$.
Poisson-primeness of $P$ then implies that some $Q_j\subseteq P$, whence
$P=Q_j$, proving that $P$ is prime. \qed\enddemo

\definition{6.3\. Poisson-primitive ideals and symplectic cores} Let $R$
be a  (commutative) Poisson algebra.

The {\it Poisson-primitive ideals\/} of $R$ are the Poisson cores of the
maximal ideals of $R$. Note from \S6.1 that all Poisson-primitive ideals
are Poisson-prime.

This terminology is chosen to reflect the following parallel.
An ideal $P$ in an algebra $A$ is left primitive if and only if $P$ is
the largest ideal contained in some maximal left ideal. If we view
$A$ as a (noncommutative) Poisson algebra via the commutator bracket $[-,-]$, then the
ideals of $A$ are precisely the Poisson left ideals. Thus, the left
primitive ideals of $A$ are exactly the Poisson cores of the maximal left
ideals. 

The {\it Poisson-primitive spectrum\/} of $R$, denoted $\Pprim R$, is the
set of all Poisson-primitive ideals of $R$. This is a subset of $\Pspec
R$, and we give it the relative topology.

By definition, the process of taking Poisson cores defines a surjective
map
$$\maxspec R \longrightarrow \Pprim R,$$
and we note that this map is continuous. Its fibres, namely the sets
$$\{ \mm\in\maxspec R \mid \P(\mm)= P\}$$
for $P\in\Pprim R$, are
called {\it symplectic cores\/}. They determine a partition of $\maxspec
R$.

Now suppose that $R=\O(V)$ is the coordinate ring of an affine variety
$V$, and that $k$ is algebraically closed. As in the complex case, we say
that
$V$ is a {\it Poisson variety\/}. Since $k$ is algebraically closed,
there is a natural identification
$V\equiv \maxspec R$, with which we transfer the symplectic cores from
$\maxspec R$ to $V$. In other words, the {\it symplectic cores in $V$\/}
are the sets
$$\{ p\in V \mid \P(\mm_p)= P\}$$
for $P\in\Pprim R$, where $\mm_p= \{f\in R \mid f(p)=0\}$.
\enddefinition 

\definition{6.4\. Example} Return to Example 3.8, and set $R= \Ogs=
\CC[x_1,x_2,x_3]$. The Poisson-primitive ideals of $R$ can be computed as
follows:
$$\xalignat2
\P(\langle x_1-\alpha,x_2,x_3\rangle) &= \langle
x_1-\alpha,x_2,x_3\rangle &&(\alpha\in\CC)  \\
\P(\langle x_1-\alpha,x_2-\beta,x_3\rangle) &= \langle x_3\rangle
&&(\alpha\in\CC,\ \beta\in\CCx)  \\
\P(\langle x_1-\alpha,x_2,x_3-\gamma\rangle) &= \langle x_2\rangle
&&(\alpha\in\CC,\ \gamma\in\CCx)  \\
\P(\langle x_1-\alpha,x_2-\beta,x_3-\gamma\rangle) &= \langle0\rangle
&&(\alpha\in\CC,\ \beta,\gamma\in\CCx) .
\endxalignat$$
It follows that the symplectic cores in $\gg^*$ are the sets
\roster
\item"$\bullet$" the individual points on the $x_1$-axis;
\item"$\bullet$" the $x_1x_2$-plane minus the $x_1$-axis;
\item"$\bullet$" the $x_1x_3$-plane minus the $x_1$-axis;
\item"$\bullet$" $\gg^*$ minus the $x_1x_2$- and $x_1x_3$-planes.
\endroster
These are precisely the coadjoint orbits of the adjoint algebraic group
of $\gg$, as we saw in Example 5.6.
\enddefinition

\head 7. Symplectic cores versus symplectic leaves \endhead

Symplectic cores are related to symplectic leaves by the following result
of Brown and Gordon \cite{\BrGr, Lemma 3.3 and Proposition 3.6}; further
relations will be given below. Here ``locally closed'' refers to the
Zariski topology.

\proclaim{7.1\. Theorem} {\rm [Brown-Gordon]} Let $V$ be a complex affine
Poisson variety.

{\rm (a)} Each symplectic core in $V$ is locally closed, and is a union
of symplectic leaves.

{\rm (b)} If the symplectic leaves in $V$ are all locally closed, then
they coincide with the symplectic cores.
\endproclaim

It is a standard result that the orbits of a connected algebraic group $G$
acting on a variety $X$ can be recovered from the orbit closures, as
follows. Take any orbit closure $C$, and remove all orbit closures
properly contained in $C$. The result will be a single $G$-orbit, and all
$G$-orbits in $X$ are obtained by this means. Yakimov has conjectured
that the symplectic cores in a complex affine Poisson variety can be
recovered from the closures of the symplectic leaves in a similar manner.
We verify this below, with the help of the following lemma of Brown and
Gordon \cite{\BrGr, Lemma 3.5}. All topological properties are to be taken
relative to the Zariski topology.

\proclaim{7.2\. Lemma} {\rm [Brown-Gordon]}  Let $V$ be a complex affine
Poisson variety, and $R=\O(V)$. Let $L$ be a symplectic leaf in $V$, and
set $K= \{f\in R \mid f=0 \text{\ on\ } L\}$. Then $K$ is a Poisson-primitive ideal of $R$, and $L$ is contained in the corresponding
symplectic core, that is, $\P(\mm_p)= K$ for all $p\in L$. \endproclaim

\proclaim{7.3\. Lemma} Let $V$ be a complex affine
Poisson variety, and $R=\O(V)$. Let $K$ be a Poisson ideal of $R$, and
$X$ the closed subvariety of $V$ determined by $K$. Then $X$ is a union
of symplectic cores and a union of symplectic leaves. In particular, the
closure of any symplectic leaf of $V$ is a union of symplectic leaves.
\endproclaim

\demo{Proof} If $p\in X$, then $\mm_p \supseteq K$. Since $K$ is a
Poisson ideal, it must be contained in the Poisson-primitive ideal $P=
\P(\mm_p)$. Now the set $C= \{q\in V\mid \P(\mm_q)= P\}$ is the
symplectic core containing
$p$, and $C\subseteq X$ because $\mm_q\supseteq P\supseteq K$ for all
$q\in C$. Therefore $X$ is a union of symplectic cores. That $X$ is a
union of symplectic leaves now follows from Theorem 7.1(a).

For any symplectic leaf $L$ of $V$, the ideal $I$ of functions in $R$ that
vanish on $L$ is a Poisson ideal by Lemma 7.2. The closed subvariety
determined by $I$ is the closure of $L$, and this is a union of symplectic
leaves by what we have just proved.
\qed\enddemo

We can now prove that symplectic cores are obtained from symplectic
leaves in the manner proposed by Yakimov; this is parts (c) and (e) of the
following theorem. Here overbars denote closures.

\proclaim{7.4\. Theorem} Let $V$ be a complex affine
Poisson variety, and $L$ a symplectic leaf in $V$.

{\rm (a)} There is a unique symplectic core $C$ in $V$ containing $L$,
and $C\subseteq \Lbar$.

{\rm (b)} $C$ is the union of those symplectic leaves of $V$ which are
dense in $\Lbar$.

{\rm (c)} $C= \Lbar \setminus \bigcup_M \Mbar$ where $M$ runs over those
symplectic leaves whose closures are properly contained in $\Lbar$.

{\rm (d)} $C$ is the unique symplectic core dense in $\Lbar$.

{\rm (e)} Each symplectic core in $V$ is dense in the closure of every
symplectic leaf it contains. Hence, it can be obtained from the closure
of such a leaf as in part {\rm (c)}.
\endproclaim

\demo{Proof} Set $R=\O(V)$, and let $K$ be the ideal of functions in $R$
that vanish on $L$.

(a) The symplectic cores and the symplectic leaves both
partition $V$, and the latter form a finer partition, by Theorem 7.1(a).
This implies the existence and uniqueness of $C$.

By Lemma 7.2, $K$ is a Poisson-primitive ideal, and the
symplectic core it determines contains $L$. By uniqueness, this core is
$C$, that is, $C= \{p\in V\mid \P(\mm_p)= K\}$. In particular,
$\mm_p\supseteq K$ for all $p\in C$, from which it follows that
$C\subseteq \Lbar$. 

(b) If $M$ is a
symplectic leaf which is dense in $\Lbar$, then $K$ equals the ideal of
functions in $R$ that vanish on $M$, and Lemma 7.2 implies that
$M \subseteq C$. On the other hand, if $M'$ is a symplectic leaf which is
contained in but not dense in $\Lbar$, the ideal $K'$ of functions
vanishing on $M'$ properly contains $K$, whence $M'$ is contained in a
symplectic core different from $C$. In this case, $M'$ is disjoint from
$C$. Part (b) now follows, because $\Lbar$ is a union of symplectic
leaves, by Lemma 7.3.

(c) In view of Lemma 7.3, the given union $\bigcup_M \Mbar$ equals the
union of those symplectic leaves which are contained in
$\Lbar$ but not dense in $\Lbar$. The given formula for $C$ thus follows
from part (b).

(d) Clearly $C$ is dense in $\Lbar$, since $L\subseteq C\subseteq \Lbar$.
If $D$ is a different symplectic core contained in $\Lbar$, then by (b),
any symplectic leaf $N\subseteq D$ is not dense in $\Lbar$. But
$D\subseteq \Nbar$ by (a), and thus $D$ is not dense in $\Lbar$.

(e) Suppose that $D$ is a symplectic core in $V$, and $N$ a symplectic
leaf contained in $D$. By (a), $D$ is the unique symplectic core
containing $N$, and $D\subseteq \Nbar$, whence $D$ is dense in $\Nbar$.
The final statement now follows from (c), with $C$ and $L$ replaced by
$D$ and $N$. \qed\enddemo

\head 8. Symplectic cores versus primitive ideals for solvable Lie
algebras \endhead

We now show that the concept of symplectic cores exactly overcomes the
limitations of symplectic leaves with respect to the Dixmier map
discussed in Section 5. Namely, the Dixmier map provides a
homeomorphism from the space of symplectic cores in $\gg^*$ onto the
primitive spectrum of $\Ug$, for any solvable finite dimensional complex
Lie algebra $\gg$. This just amounts to showing that the coadjoint orbits
in $\gg^*$, with respect to the adjoint algebraic group of $\gg$, coincide
with the symplectic cores. Solvability is not needed for the latter
result.

All that is required to obtain the new statement about the Dixmier map is
to reinterpret parts of the development of Theorem 5.4 in terms of the new
concepts. This reinterpretation also shows that (for $\gg$ solvable)
$\Pprim\Ogs$ is homeomorphic to $\prim\Ug$. With a little extra effort, we
can handle prime ideals as well, showing that
$\Pspec\Ogs$ is homeomorphic to $\spec\Ug$.

Throughout this section, $\gg$ will denote a finite dimensional complex
Lie algebra and
$G$ its adjoint algebraic group. We do not assume $\gg$ solvable until
Theorem 8.5. Some of the results we will need are developed in the
literature in terms of
$\Sg$ rather than
$\Ogs$. This requires use of the Poisson isomorphism $\theta: \Sg @>{\
\cong\ }>> \Ogs$ of (2.6).

\definition{8.1\. Actions of $G$ and $\gg$} The group $G$ acts on $\gg$
and $\gg^*$ by the adjoint and coadjoint actions, respectively, as in
\S5.1. In turn, these induce actions of $G$ by $\CC$-algebra
automorphisms on $\Sg$ and $\Ogs$, actions which we also refer to as {\it
adjoint\/} and {\it coadjoint actions\/}. All $G$-actions we mention will
refer to one of these four cases. Let us write $\Gspec \Sg$ and $\Gspec
\Ogs$ for the sets of $G$-stable prime ideals in $\Sg$ and $\Ogs$,
respectively, equipped with the relative topologies from $\spec\Sg$ and
$\spec\Ogs$.

We claim that the isomorphism $\theta$ is $G$-equivariant. To
see this, let $\{e_1,\dots,e_n\}$ be a basis for $\gg$ and
$\{\alpha_1,\dots,\alpha_n\}$ the corresponding dual basis for $\gg^*$.
As in \S2.6, $\Ogs= \CC[x_1,\dots,x_n]$ where each $x_i= \theta(e_i)$.
Given $\gamma \in G$, there are scalars $\gamma_{ij} \in \CC$ such that
$\gamma.e_j= \sum_i
\gamma_{ij}e_i$ for all $j$. Consequently,
$$(\gamma.x_j)(\alpha_i)= x_j(\gamma^{-1}.\alpha_i)=
(\gamma^{-1}.\alpha_i)(e_j)= \alpha_i(\gamma.e_j)= \gamma_{ij}$$ 
for all $i$, $j$, from which we conclude that $\gamma.x_j= \sum_i
\gamma_{ij}x_i$ for all $j$. Therefore $\gamma.\theta(e_j)=
\theta(\gamma.e_j)$ for all $j$, and the $G$-equivariance of $\theta$
follows.

For each $e\in \gg$, the Lie derivation $\ad e= [e,-]$ on $\gg$ extends
uniquely to a derivation on $\Sg$, namely the Hamiltonian $\{e,-\}$. This
yields an action of $\gg$ on $\Sg$ by derivations. We write $\gspec \Sg$
for the set of $\gg$-stable prime ideals of $\Sg$, equipped with the
relative topology from $\spec\Sg$.
\enddefinition

\proclaim{8.2\. Lemma} {\rm (a)} $\Gspec\Sg= \gspec\Sg= \Pspec\Sg$.

{\rm (b)}
$\Gspec\Ogs= \Pspec\Ogs$.

{\rm (c)} $\theta$ induces a homeomorphism
$\gspec\Sg @>{\;\approx\;}>> \Pspec\Ogs$.
\endproclaim

\demo{Proof} (a) Since $\gg$ generates the algebra $\Sg$, the
$\gg$-stable ideals of $\Sg$ coincide with the Poisson ideals. Hence,
$\gspec\Sg=
\Pspec\Sg$. By \cite{\BGR, \S13.1} or \cite{\TY, \S24.8.3}, the $\gg$-stable
ideals of $\Sg$ coincide with the $G$-stable ideals. From this,
we immediately obtain $\Gspec\Sg= \gspec\Sg$.

(b)(c) These follow immediately from (a), because $\theta$ is both
$G$-equivariant and a Poisson isomorphism. \qed\enddemo

Following our previous notation for maximal ideals corresponding to
points in varieties, write $\mm_\alpha$ for the maximal ideal of $\Ogs$
corresponding to a point $\alpha\in\gg^*$.

\proclaim{8.3\. Proposition} Let $\gg$ be a finite dimensional complex
Lie algebra and $G$ its adjoint algebraic group. There is a homeomorphism
$\phi:
\gg^*/G \rightarrow \Pprim\Ogs$ such that $\phi(G.\alpha)=
\P(\mm_\alpha)$ for all $\alpha\in\gg^*$.
\endproclaim

\demo{Proof} Since $\Sg$ is isomorphic to $\Ogs$, its maximal ideal space
is homeomorphic to $\gg^*$. A coordinate-free way to express the inverse
isomorphism is to send each $\alpha\in\gg^*$ to the ideal
$\mun_\alpha= \langle e-\alpha(e) \mid e\in\gg \rangle$ of $\Sg$. Observe
that $\theta(\mun_\alpha)= \mm_\alpha$.

By \cite{\BGR, Lemma 13.2 and proof}, there is a topological embedding
$$\tau: \gg^*/G \longrightarrow \gspec \Sg$$
such that $\tau(G.\alpha)= \bigcap_{\gamma\in G} \gamma.\mun_\alpha$ for
$\alpha\in \gg^*$. Thus, $\tau(G.\alpha)$ is the largest $G$-stable ideal
of $\Sg$ contained in $\mun_\alpha$. Invoking \cite{\BGR, \S13.1} or
\cite{\TY,
\S24.8.3} again, we find that $\tau(G.\alpha)$ is the largest
$\gg$-stable ideal of $\Sg$ contained in $\mun_\alpha$. In particular, it
now follows from \cite{\Dbook, Lemma 3.3.2} that $\tau(G.\alpha)$ is a
prime ideal. Hence, we can say that $\tau(G.\alpha)$ equals the largest
member of $\Gspec\Sg$ contained in $\mun_\alpha$.
Since $\theta$ is $G$-equivariant, it follows
that $\theta\tau(G.\alpha)$ equals the largest member of $\Gspec\Ogs$
contained in $\mm_\alpha$. In view of Lemma 8.2(b), we conclude that
$\theta\tau(G.\alpha)= \P(\mm_\alpha)$.

Combining the above with Lemma 8.2(c), we obtain a topological embedding
$$\phi: \gg^*/G \rightarrow \Pspec\Ogs$$
such that  $\phi(G.\alpha)=
\P(\mm_\alpha)$ for $\alpha\in\gg^*$. Since the image of $\phi$ is, by
definition, $\Pprim\Ogs$, the proposition is proved. \qed\enddemo

\proclaim{8.4\. Corollary} Let $\gg$ be a finite dimensional complex
Lie algebra and $G$ its adjoint algebraic group. The $G$-orbits in $\gg^*$ are precisely the
symplectic cores. \endproclaim

\demo{Proof} Injectivity and well-definedness of the homeomorphism $\phi$
of Proposition 8.3 say that for all $\alpha,\beta\in \gg^*$, we have
$G.\alpha= G.\beta$ if and only if $\P(\mm_\alpha)= \P(\mm_\beta)$. Thus,
$\alpha$ and $\beta$ lie in the same $G$-orbit if and only if they lie in
the same symplectic core. \qed\enddemo

Corollary 8.4 allows us to phrase the
Dixmier-Conze-Duflo-Rentschler-Mathieu Theorem in terms of symplectic
cores:

\proclaim{8.5\. Theorem} Let $\gg$ be a solvable finite dimensional
complex Lie algebra, and let $X$ be the set of symplectic cores in
$\gg^*$, with the quotient topology induced from $\gg^*$. Then the
Dixmier map induces a homeomorphism $X \rightarrow \prim\Ug$.
\qed\endproclaim

\definition{8.6\. The extended Dixmier map} Continue to assume that $\gg$
is solvable. Via the embedding
$\gg^*/G
\longrightarrow \gspec \Sg$ from \cite{\BGR, Lemma 13.2} used above,
identify $\gg^*/G$ with a subspace of $\gspec\Sg$. Borho, Gabriel, and
Rentschler showed that $\Dxbar$ extends uniquely to a continuous map
$$\Dxtil: \gspec \Sg \rightarrow \spec \Ug,$$
given by the rule 
$$\Dxtil(P)= \bigcap\, \{\Dx(\alpha) \mid \alpha\in\gg^*
\text{\ and\ } \mun_\alpha \supseteq P\}$$
for $P\in \gspec \Sg$
\cite{\BGR, Satz 13.4}. They named this the {\it extended Dixmier
map\/}, and proved that it is a continuous bijection \cite{\BGR, Satz
13.4, Kor\. 15.1}. Their methods, combined with Mathieu's theorem, imply
that $\Dxtil$ is a homeomorphism, as we will see shortly.
\enddefinition

\definition{8.7\. Quasi-homeomorphisms and sauber spaces} Let $X$ and $Y$
be topological spaces. A continuous map
$\phi: X\rightarrow Y$ is a {\it quasi-homeomorphism\/} provided the
induced map $F\mapsto \phi^{-1}(F)$ is an isomorphism from the lattice of
closed subsets of $Y$ onto the lattice of closed subsets of $X$. If $X$
is a subspace of $Y$, the inclusion map $X\rightarrow Y$ is a
quasi-homeo\-morphism if and only if $\overline{F\cap X}= F$ for all
closed sets $F\subseteq Y$ \cite{\BGR, \S1.6}. Borho, Gabriel, and Rentschler observed that the
inclusion map $\prim\Ug \rightarrow \spec\Ug$ is a quasi-homeomorphism
\cite{\BGR, Beispiel 1.6}, as is the above embedding $\gg^*/G
\longrightarrow \gspec \Sg$ \cite{\BGR, Lemma 13.2}.

A {\it generic point\/} of a closed subset $F\subseteq X$ is any point
$x\in F$ such that $F= \overline{\{x\}}$. The space $X$ is {\it sauber\/}
(English: {\it tidy\/}) provided every irreducible closed subset of $X$
has precisely one generic point. As observed in \cite{\BGR, \S13.3}, the
prime spectrum of any noetherian ring is sauber. We include the short
argument in the lemma below, for the reader's convenience. The same argument
shows that $\gspec \Sg$ is sauber. These spaces are noetherian as well,
since they have Zariski topologies arising from noetherian rings.
\enddefinition

\proclaim{8.8\. Lemma} Let $A$ be a noetherian ring and $R$ a
commutative noetherian Poisson $k$-al\-ge\-bra, with $\chr k=0$.

{\rm (a)} The prime spectrum $\spec A$ is a sauber noetherian space, and
if $A$ is a Jacobson ring, the inclusion map $\prim A \rightarrow \spec A$
is a quasi-homeomorphism.

{\rm (b)} The Poisson prime spectrum $\Pspec R$ is a sauber noetherian
space, and if
$R$ is an affine $k$-algebra, the inclusion map $\Pprim R\rightarrow
\Pspec R$ is a quasi-homeomorphism. \endproclaim

\demo{Proof} (a) Suppose that $F_1 \supseteq F_2 \supseteq \cdots$ is a
decreasing sequence of closed sets in $\spec A$. We may write each $F_j=
V(I_j)$ where $I_j= \bigcap F_j$. Then $I_1\subseteq I_2\subseteq \cdots$
is an increasing sequence of ideals of $A$. Since this sequence
stabilizes, so does the original sequence of closed sets. Thus, $\spec A$
is a noetherian space.

Let $F= V(I)$ be an arbitrary closed subset of $\spec A$,
where $I$ is an ideal of $A$. We may replace $I$ by its prime radical,
so there is no loss of generality in assuming that $I$ is semiprime. Since
$A$ is noetherian, there are only finitely many prime ideals minimal
over $I$, say $Q_1,\dots,Q_n$, and $I= Q_1\cap \cdots\cap Q_n$. It follows
that $F= V(Q_1)\cup \cdots\cup V(Q_n)$. 

If $F$ is irreducible, then $F= V(Q_j)$ for some $j$. In this case, $Q_j$
is the unique generic point of $F$, proving that $\spec A$ is sauber.

Now assume that $A$ is a Jacobson ring, so that all prime ideals of $A$
are intersections of primitive ideals. It follows that
$$I= \bigcap\, F= \bigcap\, (F\cap \prim A),$$
from which we see that $F$ equals the closure of $F\cap \prim A$ in
$\spec A$. Thus, by \cite{\BGR, \S1.6}, the inclusion map $\prim A
\rightarrow \spec A$ is a quasi-homeomorphism. 

(b) The argument applied in (a) also shows that $\Pspec R$ is a
noetherian space.

As discussed in \S6.1, any closed set $F$ in $\Pspec R$ can be
written $F= V_P(I)$ for some Poisson ideal $I$. There are only finitely
many prime ideals minimal over $I$, say $Q_1,\dots,Q_n$, and the $Q_i$
are Poisson ideals by Lemma 6.2. Hence, we may replace $I$ by $Q_1\cap
\cdots\cap Q_n$, and it follows that $F= V_P(Q_1)\cup \cdots\cup
V_P(Q_n)$.

Just as in (a), if $F$ is irreducible, $F= V_P(Q_j)$ for some $j$, and
then $Q_j$ is the unique generic point of $F$. This proves that $\Pspec
R$ is sauber.

Now assume that $R$ is an affine $k$-algebra. Then $R$ is a Jacobson
ring, and it follows that every Poisson prime ideal of $R$ is an
intersection of Poisson-primitive ideals (e.g., see \cite{\GAth, Lemma
1.1(e)}). From this, we conclude as in (a) that the inclusion map $\Pprim
R\rightarrow
\Pspec R$ is a quasi-homeomorphism. 
\qed\enddemo

\proclaim{8.9\. Lemma} Let $X\subseteq X'$ and $Y\subseteq Y'$ be
topological spaces, such that $X'$ and $Y'$ are sauber and noetherian.
Assume also that the inclusion maps $X\rightarrow X'$ and $Y\rightarrow
Y'$ are quasi-homeomorphisms. Then any continuous map $\phi: X\rightarrow
Y$ extends uniquely to a continuous map $\phi': X'\rightarrow Y'$.
Moreover, if $\phi$ is a homeomorphism, so is $\phi'$. \endproclaim

\demo{Proof} The existence and uniqueness of $\phi$ are proved in
\cite{\BGR, Lemma 13.3}. The final statement follows by the usual
universal property argument. \qed\enddemo

\proclaim{8.10\. Theorem} {\rm [Borho-Gabriel-Rentschler-Mathieu]} Let
$\gg$ be a solvable finite dimensional complex Lie algebra. The
extended Dixmier map
$$\Dxtil: \gspec \Sg \longrightarrow \spec \Ug$$
is a homeomorphism. \endproclaim

\demo{Proof} Following the proof of \cite{\BGR, Satz 13.4}, recall that
$\gspec \Sg$ and $\spec \Ug$ are sauber noetherian spaces, and that the
embedding $\gg^*/G \rightarrow \gspec \Sg$ and the inclusion $\prim \Ug
\rightarrow \spec \Ug$ are quasi-homeomorphisms. The map $\Dxtil$ is
defined, with the help of Lemma 8.9, to be the unique continuous map
from $\gspec \Sg$ to $\spec \Ug$ extending $\Dxbar$. Since $\Dxbar$ is a
homeomorphism, Lemma 8.9 implies that $\Dxtil$ is a homeomorphism.
\qed\enddemo

In Poisson-ideal-theoretic terms, Theorems 5.4 and 8.10 can be restated as
follows.

\proclaim{8.11\. Theorem}  Let $\gg$ be a solvable finite dimensional
complex Lie algebra. Then  there is a homeomorphism 
$$\psi: \Pprim\Ogs \longrightarrow \prim\Ug$$
such  that $\psi(\P(\mm_\alpha))= \Dx(\alpha)$ for
$\alpha\in\gg^*$, and $\psi$ extends uniquely to a homeomorphism 
$$\Pspec \Ogs \longrightarrow \spec \Ug.$$
\endproclaim

\demo{Proof} To obtain $\psi$, just compose the factorized Dixmier map
$\Dxbar$ with the inverse of the homeomorphism $\phi$ of Proposition 8.3.
By Lemma 8.8, $\Pspec\Ogs$ and $\spec\Ug$ are sauber noetherian spaces,
and the inclusion maps $\Pprim\Ogs
\rightarrow \Pspec\Ogs$ and $\prim\Ug \rightarrow
\spec\Ug$ are quasi-homeomorphisms. Therefore the existence and uniqueness
of the desired extension of $\psi$ follow from Lemma 8.9. \qed\enddemo

\head 9. Modified conjectures for quantized coordinate rings \endhead

In light of Theorems 7.1, 7.4, 8.5, and 8.11, we nominate the concept of
symplectic cores as the {\it best algebraic approximation\/} for
symplectic leaves. Further, we suggest that symplectic leaves should be
replaced by symplectic cores in applications of the Orbit Method to
algebraic problems. In particular, we revise and refine the general
principle discussed in \S4.4 to the following conjecture. It is, of
necessity, somewhat imprecise, given the lack of a precise definition of
the concept of quantized coordinate rings.

\definition{9.1\. Primitive spectrum conjecture for quantized coordinate
rings} 
\roster
\item"" {\sl Assume that $k$ is algebraically closed of characteristic
zero, and let $A$ be a generic quantized coordinate ring of an affine
algebraic variety $V$ over $k$. Then $A$ should be a member of a flat
family of $k$-algebras with semiclassical limit $\O(V)$, such that $\prim
A$ is homeomorphic to the space of symplectic
cores in $V$, with respect to the semiclassical limit Poisson structure.
Further, there should be compatible homeomorphisms $\prim A \rightarrow
\Pprim \O(V)$ and $\spec A \rightarrow \Pspec \O(V)$.}
\endroster
Each of the known types of quantized coordinate rings supports an
action of an algebraic torus $H= (\kx)^m$ (see \cite{\BrGd, \S\S
II.1.14-18} for a summary), which has a parallel action (by Poisson
automorphisms) on the semiclassical limit (e.g., see
\cite{\GoYa, \S0.2; \GLau, Section 2}). We tighten the conjecture
above and posit that there should exist homeomorphisms as desribed which
are also equivariant with respect to the relevant torus actions.
\enddefinition

\definition{9.2\. Remarks} {\bf (a)} The discussion of the simple
example $A_q= \Oq(\CC^2)$ in \S4.5 indicates why Conjecture 9.1 is
restricted to generic quantized coordinate rings. In particular, $\prim
A_q$ has a generic point when $q$ is not a root of unity, but no generic
points otherwise. Since $\Pspec \CC[x,y]$ has a generic point, it is not
homeomorphic to $\prim A_q$ when $q$ is a root of unity.

{\bf (b)} Each of the ``standard'' single parameter quantized coordinate
rings is defined as a member of a one-parameter family of algebras, and
it is this (flat) family to which the conjecture is meant to apply. For
instance, the algebras $\Oq(SL_n(k))$ (with $n$ fixed) are defined for
all $q\in\kx$ in the same way (e.g., \cite{\BrGd, \S I.2.4}), and
substituting an indeterminate
$t$ for
$q$ in the definition results in a torsionfree $k[t^{\pm1}]$-algebra $A$
with $A/(t-q)A \cong \Oq(SL_n(k))$ for all $q\in\kx$, just as with the
case $n=2$ in \S\S 1.6, 2.2(c). The semiclassical limit is $\O(SL_n(k))$
with the Poisson bracket satisfying
$$\xalignat4
\{X_{ij},X_{im}\} &= X_{ij}X_{im} &&(j<m)  \\
\{X_{ij},X_{lj}\} &= X_{ij}X_{lj} &&(i<l)  \\
\{X_{ij},X_{lm}\} &= 0 &&(i<l,\ j>m)  \\
\{X_{ij},X_{lm}\} &= 2X_{im}X_{lj} &&(i<l,\ j<m).
\endxalignat$$
This Poisson structure and the above flat family should feature in the
$SL_n$ case of Conjecture 9.1, that is, for $q$ not a root of unity,
$\prim \Oq(SL_n(k))$ should be homeomorphic to the space of symplectic
cores in $SL_n(k)$ and to $\Pprim \O(SL_n(k))$, and $\spec \Oq(SL_n(k))$
should be homeomorphic to $\Pspec \O(SL_n(k))$. Such a ``standard''
version of the conjecture is to be posed for
$\Oq(M_n(k))$,
$\Oq(GL_n(k))$,
$\Oq(G)$, and other ``standard'' cases.

The situation is more involved for ``nonstandard'' cases, and for
multiparameter families, which have to be reduced to single parameter
families in order to obtain semiclassical limits. In such cases, the
conjecture may be sensitive to the choice of flat family -- different
flat families may yield different Poisson structures in the semiclassical
limit, and the conjecture may hold for some of these semiclassical limits
but not for others. This phenomenon appears in an example of Vancliff
\cite{\Vnc, Example 3.14}, which we discuss in Example 9.9.

{\bf (c)} As discussed at the end of Example 2.6, the enveloping
algebra of a finite dimensional Lie algebra $\gg$ is a generic member of the flat family given by the
homogenization of $\Ug$, and so $\Ug$ should qualify as a
generic quantized coordinate ring of
$\gg^*$. The semiclassical limit of this family is the
Poisson algebra $\Ogs$. For this setting, K\. A\. Brown has noted
difficulties with Conjecture 9.1 in what one might expect to be the most
canonical case, namely when $\gg$ is semisimple
\cite{\Bro}. Following the Orbit Method, one would seek a
bijection $\L \longleftrightarrow P$ between symplectic leaves in $\gg^*$
and primitive ideals in $\Ug$ such that the Gelfand-Kirillov dimension of
$\Ug/P$ equals the dimension of $\L$. In particular, the
zero-di\-men\-sion\-al symplectic leaves of $\gg^*$, which are the same as
the zero-dimensional symplectic cores, should match up with the maximal
ideals of finite codimension in
$\Ug$. However, $\Ug$ has infinitely many such maximal ideals, while
there is only one zero-dimensional symplectic leaf in $\gg^*$. (The
latter can be verified by using Theorem 4.2 together with the fact that
the identification of $\gg^*$ with $\gg$ via the Killing form identifies
the coadjoint orbits in $\gg^*$ with the adjoint orbits in $\gg$
\cite{\CoMc, p\. 12}.)

Other differences are already visible in the case $\gg=\frak{sl}_2(\CC)$.
As is easily computed, all but one of the coadjoint orbits in $\gg^*$ are
closed (compare with the adjoint orbits, computed in \cite{\CoMc, Example
1.2.1}). It follows (using Proposition 8.3, or by direct computation) that
all but one of the points of $\Pprim
\Ogs$ are closed. However, $\prim\Ug$ has infinitely many non-closed
points, and therefore it is not homeomorphic to $\Pprim \Ogs$.

{\bf (d)} Whenever Conjecture 9.1 does hold, the space of
symplectic cores in $V$ must be homeomorphic to $\Pprim \O(V)$. It is an
open question whether the space of symplectic cores for an arbitrary
affine Poisson algebra $R$ is homeomorphic to $\Pprim R$, but this does
hold when $R$ satisfies the {\it Poisson Dixmier-Moeglin equivalence\/},
as follows from \cite{\GAth, Theorem 1.5}; we excerpt the basic
argument in Lemma 9.3. This equivalence requires that the
Poisson-primitive ideals of $R$ coincide with the locally closed points
of $\Pspec R$, and with those Poisson prime ideals $P$ of $R$ for which
the Poisson center (cf\. \S9.6(b)) of the quotient field of $R/P$ is algebraic over $k$.
It holds for the semiclassical
limits of many quantized coordinate rings via
\cite{\GAth, Theorem 4.1}, as shown in \cite{\GLau, Section 2}.

{\bf (e)} As in Theorem 8.11, the existence of a homeomorphism $\prim A
\rightarrow \Pprim \O(V)$ as in the conjecture typically implies the
existence of a compatible homeomorphism $\spec A \rightarrow \Pspec
\O(V)$. We display this in Lemma 9.4 below for emphasis. On the other
hand, a homeomorphism $\spec A \rightarrow \Pspec \O(V)$ will restrict to
a homeomorphism  $\prim A
\rightarrow \Pprim \O(V)$ provided $\O(V)$ satisfies the Poisson
Dixmier-Moeglin equivalence and $A$ satisfies the {\it Dixmier-Moeglin
equivalence\/}. The latter equivalence requires that the primitive ideals
of $A$ coincide with the locally closed points of $\spec A$, and with
those prime ideals $P$ of $A$ for which $Z(\Fract A/P)$ is algebraic over
$k$. It was verified for many quantized coordinate rings in
\cite{\GLDM} (see
\cite{\BrGd, Corollary II.8.5} for a summary).
\enddefinition

\proclaim{9.3\. Lemma} Let $R$ be a commutative affine Poisson
$k$-algebra, and assume that all Poisson-primitive ideals of $R$ are
locally closed points in $\Pspec R$. Then the Zariski topology on
$\Pprim R$ coincides with the quotient topology induced by the Poisson
core map
$\P(-): \maxspec R \rightarrow \Pprim R$. Consequently, the space of
symplectic cores in $\maxspec R$ is homeomorphic to $\Pprim R$.
\endproclaim

\demo{Proof} Observe first that the map $\P(-)$ is continuous. It is
surjective by definition of $\Pprim R$. 

We claim that $P= \bigcap\, \{\mm\in\maxspec R\mid \P(\mm)=P\}$, for any
Poisson-primitive ideal $P$ of $R$. Since
$P$ is locally closed in
$\Pspec R$ (by assumption), the singleton $\{P\}$ is open in its closure
$V_P(P)$, and so $\{P\}= V_P(P) \setminus V_P(J)$ for some Poisson ideal
$J$ of $R$. Note that $J\not\subseteq P$; hence, after replacing $J$ by
$J+P$, we may assume that $J\supsetneq P$. If $\mm\supseteq P$ is a
maximal ideal such that $\P(\mm) \ne P$, then
$\mm\supseteq \P(\mm)\supseteq J$. The remaining maximal ideals
containing $P$ must intersect to $P$ by the Nullstellensatz, verifying
the claim.

Now consider a set $X\subseteq \Pprim R$ whose inverse image under
$\P(-)$, call it $Y$, is closed in $\maxspec R$. Thus,
$$Y= \{\mm\in\maxspec R\mid \P(\mm)\in X\} = \{\mm\in\maxspec R\mid \mm
\supseteq I\}$$
for some ideal $I$ of $R$. If $P\in X$ and $\mm\in\maxspec R$ with
$\P(\mm)= P$, then $\mm\in Y$, and so $\mm\supseteq I$. By the claim
above, the intersection of these maximal ideals equals $P$, and thus
$P\supseteq I$. Conversely, if $P\in\Pprim R$ and $P\supseteq I$, then
$P=\P(\mm)$ for some maximal ideal $\mm\supseteq I$, whence $\mm\in Y$
and $P\in X$. Therefore $X= \{P\in\Pprim R\mid P\supseteq I\}$, a closed
set in $\Pprim R$. This proves that the topology on $\Pprim R$ is the
quotient topology inherited from $\maxspec R$ via $\P(-)$.

The final statement of the lemma follows directly. \qed\enddemo

\proclaim{9.4\. Lemma} Let $A$ be a noetherian $k$-algebra and
$R$ a commutative noetherian Poisson $k$-algebra, with $\chr k =0$.

{\rm (a)} A bijection $\phi: \spec A \rightarrow \Pspec R$ is a
homeomorphism if and only if $\phi$ and $\phi^{-1}$ preserve inclusions.

{\rm (b)} Assume that $A$ is a Jacobson ring and $R$ an affine
$k$-algebra. Then any homeomorphism
$\prim A\rightarrow \Pprim R$ extends uniquely to a homeomorphism $\spec A
\rightarrow \Pspec R$.

{\rm (c)} Assume that $A$ satisfies the Dixmier-Moeglin equivalence and
$R$ the Poisson Dixmier-Moeglin equivalence. Then any
homeomorphism $\spec A \rightarrow \Pspec R$ restricts to a homeomorphism
$\prim A\rightarrow \Pprim R$.
\endproclaim

\demo{Proof} (a) For $P,Q\in \spec A$, we have $P\subseteq Q$ if and only
if $Q\in \overline{\{P\}}$, and similarly in $\Pspec R$. Hence, any
homeomorphism between these spaces must preserve inclusions.

Conversely, if $\phi$ and $\phi^{-1}$ preserve inclusions, then
$\phi(V(P))= V_P(\phi(P))$ for all $P\in \spec A$. Since the closed sets
in $\spec A$ are exactly the finite unions of $V(P)$\,s (recall the proof
of Lemma 8.8(a)), it follows that $\phi$ sends closed sets to closed
sets, i.e., $\phi^{-1}$ is continuous. Similarly, $\phi$ is continuous,
and hence a homeomorphism.

(b) Lemmas 8.8 and 8.9.

(c) Under the assumed equivalences, $\prim A$ consists of the locally
closed points in $\spec A$, and $\Pprim R$ consists of the locally closed
points in $\Pspec R$. \qed\enddemo

\definition{9.5\. Example} Let $A_q= \Oq(k^2)$, where $k=\kbar$,
$\chr k=0$, and $q\in\kx$. View $R= \O(k^2)$ as the semiclassical
limit of the family $(A_q)_{q\in\kx}$, with the Poisson structure
exhibited in Example 2.2(a). The torus $H=(\kx)^2$ acts on $A_q$ via
algebra automorphisms and on $R$ via Poisson automorphisms so that (in
both cases) $(\alpha_1,\alpha_2).x_i= \alpha_ix_i$ for
$(\alpha_1,\alpha_2) \in H$ and $i=1,2$.

Assume that $q$ is not a root of unity. As is easily checked (e.g.,
\cite{\BrGd, Example II.1.2}), the prime ideals of $A_q$ are
\roster
\item"$\bullet$" the maximal ideals $\langle x_1-\alpha,\, x_2\rangle$ and
$\langle x_1,\, x_2-\beta\rangle$, for $\alpha,\beta\in k$;
\item"$(\diamondsuit)\ \bullet$" the height $1$ primes $\langle
x_1\rangle$ and
$\langle x_2\rangle$;
\item"$\bullet$" the zero ideal.
\endroster
All of these prime ideals, except for $\langle x_1\rangle$ and $\langle
x_2\rangle$, are primitive \cite{\BrGd, Example II.7.2}. The closed sets
in $\spec A_q$ are easily found, but we shall not list them here -- see
\cite{\BrGd, Example II.1.2 and Exercise II.1.C}. 

With very similar computations, one finds the Poisson prime and Poisson-primitive ideals in $R$, and a list of the closed subsets of $\Pspec R$.  
In terms of notation, the answers are the same as for $A_q$ -- the list
$(\diamondsuit)$ also describes the Poisson prime ideals of $R$, and all
of these ideals, except for $\langle x_1\rangle$ and $\langle
x_2\rangle$, are Poisson-primitive. We conclude that there exist
compatible homeomorphisms $\prim A_q \rightarrow
\Pprim R$ and $\spec A_q \rightarrow \Pspec R$, sending the ideal
$\langle x_1-\alpha,\, x_2\rangle$ of $A_q$ to the ideal $\langle
x_1-\alpha,\, x_2\rangle$  of $R$, and so on. (We say that these maps are
given by ``preservation of notation''.) These homeomorphisms are
equivariant with respect to the actions of
$H$ described above.

By inspection, all Poisson-primitive ideals of $R$ are locally closed in
$\Pspec R$. Consequently, we conclude from Lemma 9.3 that the space of
symplectic cores in
$\maxspec R \approx k^2$ is homeomorphic to $\Pprim R$.
\enddefinition

Analyzing the prime ideals in a quantized coordinate ring typically
involves investigating localizations of factor algebras, which often turn
out to be quantum tori. We sketch some basic procedures used to determine
prime ideals in quantum tori, and similar ones for the analogous ``Poisson
tori''.

\definition{9.6\. Some computational tools} {\bf (a)} A {\it quantum
torus\/} over $k$ is the localization of a quantum affine space
$\O_{\bfq}(k^n)$ obtained by inverting the generators $x_i$, that is, an
algebra
$$\O_{\bfq}((\kx)^n)= k\langle x^{\pm1}_1,\dots,x^{\pm1}_n \mid x_ix_j=
q_{ij}x_jx_i \text{\ for all\ } i,j \rangle,$$
where $\bfq= (q_{ij})$ is a multiplicatively antisymmetric $n\times n$
matrix over $k$. Set $T= \O_{\bfq}((\kx)^n)$.

Since $T$ is a $\ZZ^n$-graded algebra, with $1$-dimensional homogeneous components,
its center is spanned by central monomials \cite{\Hod, Lemma 1.1}. The
latter are easily computed: a monomial $x_1^{m_1} x_2^{m_2} \cdots
x_n^{m_n}$ is central if and only if $\prod_{j=1}^n q_{ij}^{m_j} =1$ for
all $i$. All ideals of $T$ are induced from ideals of $Z(T)$ \cite{\Hod,
Theorem 1.2; \qaff, Proposition 1.4}, from which it follows that
contraction and extension give inverse homeomorphisms between $\spec T$
and $\spec Z(T)$ \cite{\qaff, Corollary 1.5(b)}. In particular, it
follows from the above facts that $T$ is a simple algebra if and only if
$Z(T)=k$ \cite{\McPe, Proposition 1.3}.

{\bf (b)} Let $R= k[x^{\pm1}_1,\dots,x^{\pm1}_n]$ be a Laurent polynomial
ring, equipped with a Poisson bracket such that $\{x_i,x_j\}=
\pi_{ij}x_ix_j$ for all $i$, $j$, where $(\pi_{ij})$ is an antisymmetric
$n\times n$ matrix over $k$. The results of part (a) all have Poisson
analogs for $R$, as follows.

The {\it Poisson center\/} of $R$, denoted $Z_P(R)$, is the subalgebra
consisting of those $r\in R$ for which the derivation $\{r,-\}$ vanishes.
Since the Poisson bracket on $R$ respects the $\ZZ^n$-grading, $Z_P(R)$
is spanned by the monomials it contains \cite{\Hod, Lemma 2.1; \Vnc,
Lemma 1.2(a)}. A monomial $x_1^{m_1} x_2^{m_2} \cdots
x_n^{m_n}$ is Poisson central if and only if $\sum_{j=1}^n \pi_{ij}m_j =0$
for all $i$. All Poisson ideals of
$R$ are induced from ideals of
$Z_P(R)$ \cite{\Hod, Theorem 2.2; \Vnc, Lemma 1.2(b)}, from which it
follows that contraction and extension give inverse homeomorphisms
between $\Pspec R$ and $\spec Z_P(R)$. In
particular, it follows from the above facts that $R$ is Poisson-simple
(meaning that it has no proper nonzero Poisson ideals) if and only if
$Z_P(R)=k$.
\enddefinition

\definition{9.7\. Example} Let $A_q= \Oq(SL_2(k))$, where $k=\kbar$,
$\chr k=0$, and $q\in\kx$. View $R= \O(SL_2(k))$ as the semiclassical
limit of the family $(A_q)_{q\in\kx}$, with the Poisson structure
exhibited in Example 2.2(c). The torus $H= (\kx)^2$ again acts on $A_q$
and $R$, this time so that
$$\xalignat2 (\alpha,\beta).X_{11} &= \alpha\beta X_{11} 
&(\alpha,\beta).X_{12} &= \alpha\beta^{-1} X_{12}  \\
(\alpha,\beta).X_{21} &= \alpha^{-1}\beta X_{21}  &(\alpha,\beta).X_{22}
&= \alpha^{-1}\beta^{-1} X_{22}  
\endxalignat$$
for $(\alpha,\beta) \in H$.

Now restrict $q$ to a non-root of unity. The prime ideals of $A_q$ can be
computed with the tools of \S9.6(a), as outlined in \cite{\BrGd, Exercise
II.1.D}. For instance, one checks that $A_q$ has a localization
$$A_q[X_{11}^{-1},X_{12}^{-1},X_{21}^{-1}] \cong k\langle
x^{\pm1},y^{\pm1},z^{\pm1} \mid xy=qyx,\ xz=qzx,\ yz=zy \rangle,$$
and that the center of the latter algebra is $k[(yz^{-1})^{\pm1}]$. It
follows that the prime ideals of $A_q$ not containing $X_{12}$ or
$X_{21}$ consist of $\langle0\rangle$ and $\langle X_{12}-\lambda
X_{21}\rangle$, for $\lambda\in \kx$. The full list of prime ideals of
$A_q$ is as follows:
\roster
\item"$\bullet$" the maximal ideals $\langle X_{11}-\lambda,\, X_{12},\,
X_{21},\, X_{22}-\lambda^{-1} \rangle$, for $\lambda\in\kx$;
\item"$(\spadesuit)\ \bullet$" the ideal $\langle X_{12},\, X_{21}
\rangle$;
\item"$\bullet$" the height $1$ primes $\langle X_{21}\rangle$ and
$\langle X_{12}-\lambda X_{21}\rangle$, for $\lambda\in k$;
\item"$\bullet$" the zero ideal.
\endroster
A diagram of $\spec A_q$, with inclusions marked, is given in
\cite{\BrGd, Diagram II.1.3}.

 A similar computation, using \S9.6(b), yields the Poisson
prime ideals of
$R$, which can be described exactly as in $(\spadesuit)$. This provides a
natural $H$-equivariant bijection $\phi: \spec A_q \rightarrow \Pspec R$,
given by ``preservation of notation". By inspection, $\phi$ and
$\phi^{-1}$ preserve inclusions, and thus, by Lemma 9.4(a), $\phi$ is a
homeomorphism.

The algebra $A_q$ satisfies the Dixmier-Moeglin equivalence by
\cite{\BrGd, Corollary II.8.5}, and $R$ satisfies the Poisson
Dixmier-Moeglin equivalence by \cite{\GAth, Theorem 4.3}. Therefore Lemma
9.4(c) implies that $\phi$ restricts to a homeomorphism $\prim A_q
\rightarrow \Pprim R$. In $A_q$, all prime ideals are primitive except
for $\langle X_{12},\, X_{21} \rangle$ and $\langle0\rangle$ (cf\.
\cite{\BrGd, Example II.8.6}). Similarly, in $R$ all Poisson prime ideals
are Poisson-primitive except for $\langle X_{12},\, X_{21} \rangle$ and
$\langle0\rangle$. As in the previous example, we can use Lemma 9.3 to
see that the space of symplectic cores in $\maxspec R \approx SL_2(k)$ is
homeomorphic to $\Pprim R$.
\enddefinition

\definition{9.8\. Evidence for Conjecture 9.1} In most of the
instances discussed below, $k$ is assumed to be algebraically closed of
characteristic zero.

{\bf (a)} Examples 9.5
and 9.7 are the most basic instances in which the conjecture has been
verified. In the same way (although with somewhat more effort), one can
verify it for $\Oq(GL_2(k))$. In particular, most of the work required to
determine the prime ideals in the generic $\Oq(GL_2(k))$ is done in
\cite{\BrGd, Example II.8.7}.

{\bf (b)} We next turn to the quantized coordinate rings
$\Oq(G)$ and $\O_{q,p}(G)$ over $k=\CC$, where $G$ is a connected
semisimple complex Lie group, $q\in\kx$ is not a root of unity, and $p$
is an antisymmetric bicharacter on the weight lattice of $G$ (as in
\cite{\HLT, \S3.4}). 

The Poisson structure
on $\O(G)$ resulting from the semiclassical limit process gives $G$ the
combined structure of a {\it Poisson-Lie group\/} (e.g., see
\cite{\KoSo, Chapter 1} for the concept, and \cite{\HoLeA, \S A.1} for
the result). There is a known recipe for the symplectic leaves in
$G$ in case the Poisson structure arises from the standard
quantization \cite{\HoLeA, Appendix A}, and similarly in the
multiparameter ``algebraic'' case \cite{\HLT, Theorem 1.8}. In both these
cases, it follows that the symplectic leaves are Zariski locally closed
(see \cite{\BGY, Theorem 1.9} for a more explicit statement). Hence, the
symplectic leaves in $G$ coincide with the symplectic cores (Theorem
7.1(b)). 

As discussed in \S4.4, Hodges and Levasseur put forward the conjecture
that there should be a bijection between $\prim \Oq(G)$ and the set of
symplectic leaves in $G$
\cite{\HoLeA, \S2.8, Conjecture 1}. They developed such bijections for
$G=SL_n(\CC)$ in
\cite{\HoLeB}, and for general $G$ in their work with Toro \cite{\HLT}.
More generally, Hodges, Levasseur, and Toro established a bijection
between $\prim \O_{q,p}(G)$ and the set of
symplectic leaves in $G$ in the algebraic case. All these bijections
are equivariant with respect to natural actions of a maximal torus of $G$.

Except for the case $G= SL_2(\CC)$ covered in Example 9.7, the
topological properties of the above bijections are not known. Even when
$G= SL_3(\CC)$, it is not known whether
$\prim \Oq(G)$ is homeomorphic to the space of symplectic leaves (=
cores) in $G$.

{\bf (c)} The prime and primitive spectra of general multiparameter
quantum affine spaces $\Obfq(k^n)$ were analyzed by Goodearl and Letzter
in \cite{\GLquo}, assuming $k=\kbar$ together with a minor technical
assumption (that either
$\chr k = 2$, or $-1$ is not in the subgroup $\langle q_{ij}\rangle
\subseteq \kx$). They proved that there are compatible topological
quotient maps $k^n\approx \maxspec \O(k^n) \rightarrow \prim
\Obfq(k^n)$ and $\spec \O(k^n) \rightarrow \spec
\Obfq(k^n)$, equivariant with respect to natural actions of the torus
$(\kx)^n$ \cite{\GLquo, Theorem 4.11}. Similar results were proved not
only for quantum tori $\Obfq((\kx)^n)$ \cite{\GLquo, Theorem 3.11} but
also for quantum affine toric varieties \cite{\GLquo, Theorem 6.3}. 

Oh, Park, and Shin converted these topological quotient results into the
following (assuming $\chr k=0$ and $-1 \notin \langle q_{ij}\rangle
\subseteq \kx$): For each  $\Obfq(k^n)$, there is a Poisson structure
on $\O(k^n)$ such that there are compatible homeomorphisms $\Pprim \O(k^n)
\rightarrow \prim \Obfq(k^n)$ and $\Pspec \O(k^n) \rightarrow \spec
\Obfq(k^n)$ \cite{\OPS, Theorem 3.5}. Goodearl and Letzter, finally,
showed that such homeomorphisms could be obtained for semiclassical limit
Poisson structures \cite{\GLet, Theorem 3.6}, and extended the results to
quantum affine toric varieties \cite{\GLet, Theorem 5.2}.

All these Poisson algebra structures on $\O(k^n)$ satisfy the Poisson
Dixmier-Moeglin equivalence \cite{\GAth, Example 4.6}. Hence, the space
of symplectic cores in $k^n$ is homeomorphic to
$\prim \O_{\bfq}(k^n)$, via Lemma 9.3. The symplectic cores in $k^n$ are
algebraic, whereas this does not always hold for the symplectic leaves,
as shown by Vancliff \cite{\Vnc, Corollary 3.4}. An explicit example is
computed in \cite{\GLet, Example 3.10}.

{\bf (d)} The prime and primitive spectra of the algebras
$K^{P,Q}_{n,\Gamma}(k)$ introduced by Horton
\cite{\Hrt} were analyzed by Oh in \cite{\Ohthree}. These algebras are
multiparameter quantizations of $\O(k^{2n})$, and include quantum
symplectic spaces $\Oq({\frak{sp}}\,k^{2n})$, even-dimensional quantum
euclidean spaces $\Oq({\frak{o}}\,k^{2n})$, and quantum Heisenberg
spaces, among others. Oh introduced Poisson algebra structures
$A^{P,Q}_{n,\Gamma}(k)$ on $\O(k^{2n})$, and constructed compatible
homeomorphisms $\Pprim A^{P,Q}_{n,\Gamma}(k) \rightarrow \prim
K^{P,Q}_{n,\Gamma}(k)$ and $\Pspec A^{P,Q}_{n,\Gamma}(k) \rightarrow
\spec K^{P,Q}_{n,\Gamma}(k)$, assuming the parameters involved in $P$,
$Q$, $\Gamma$ are suitably generic
\cite{\Ohthree, Theorem 4.14}.
\enddefinition

As stated in Remark 9.2(b), a quantized coordinate ring may belong to
some flat families for which Conjecture 9.1 holds and also to others for
which it fails. We outline Vancliff's example \cite{\Vnc, Example 3.14} 
illustrating this phenomenon.

\definition{9.9\. Example} {\bf (a)} Let $a_i=i-1$ for
$i=1,2,3$, and set 
$$\align
R_0 &= \CC[h][(1+a_ih)^{-1}\mid i=1,2,3]  \\
A &= R_0\langle x_1,x_2,x_3 \mid x_ix_j= r_{ij}x_jx_i \text{\
for\ } i,j=1,2,3\rangle,
\endalign$$ 
where $r_{ij}= (1+a_ih)(1+a_jh)^{-1}$ for all $i$,
$j$. This defines a flat family of $\CC$-algebras, whose semiclassical
limit is the polynomial ring $R= \CC[x_1,x_2,x_3]$ with the Poisson
bracket satisfying
$$\xalignat3 
\{x_1,x_2\} &= -x_1x_2  &\{x_1,x_3\} &= -2x_1x_3  &\{x_2,x_3\} &= -x_2x_3
\,.
\endxalignat$$
It follows from \cite{\Vnc, Corollary 3.4} that the symplectic leaves in
$\CC^3$ for this Poisson structure are algebraic; hence, they coincide
with the symplectic cores (Theorem 7.1(b)). By \cite{\GAth, Example 4.6},
$R$ satisfies the Poisson Dixmier-Moeglin equivalence, and so Lemma 9.3
implies that the space of symplectic leaves in $\CC^3$ is homeomorphic to
$\Pprim R$. 

The Poisson-primitive ideals of $R$ are listed in \cite{\Vnc, Example
3.14} (where they are labelled ``maximal Poisson ideals''). They consist
of
\roster
\item"$(\diamondsuit)\ \bullet$" the maximal ideals $\langle
x_1-\alpha,x_2,x_3\rangle$,
$\langle x_1,x_2-\beta,x_3\rangle$, $\langle x_1,x_2,x_3-\gamma\rangle$,
for $\alpha,\beta,\gamma\in \CC$;
\item"$\bullet$" the height $1$ primes $\langle x_1\rangle$, $\langle
x_2\rangle$, $\langle x_3\rangle$, and $\langle x_1x_3-\lambda
x_2^2\rangle$, for $\lambda\in\CCx$.
\endroster
We can compute them by using \S9.6(b) to find the Poisson prime ideals of
$R$ and then applying the Poisson Dixmier-Moeglin equivalence. For
instance, the Poisson center of the localization
$\CC[x^{\pm1}_1,x^{\pm1}_2,x^{\pm1}_3]$ is
$\CC[(x_1x_2^{-2}x_3)^{\pm1}]$, from which it follows that any nonzero
Poisson prime ideal of $R$ must contain either one of the $x_i$ or else
$x_1x_3-\lambda x_2^2$ for some $\lambda\in\CCx$. The full list of
Poisson prime ideals of $R$ is
\roster
\item"$\bullet$" $\langle x_1-\alpha,x_2,x_3\rangle$,
$\langle x_1,x_2-\beta,x_3\rangle$, $\langle x_1,x_2,x_3-\gamma\rangle$,
for $\alpha,\beta,\gamma\in \CC$;
\item"$\bullet$" $\langle x_1,x_2\rangle$, $\langle x_1,x_3\rangle$,
$\langle x_2,x_3\rangle$;
\item"$\bullet$" $\langle x_1\rangle$, $\langle
x_2\rangle$, $\langle x_3\rangle$, $\langle x_1x_3-\lambda
x_2^2\rangle$, for $\lambda\in\CCx$;
\item"$\bullet$" $\langle 0\rangle$.
\endroster
Inspection immediately shows that the locally closed points of $\Pspec R$
are the Poisson prime ideals listed in $(\diamondsuit)$.

{\bf (b)} A generic member of the flat family given by $A$ is $B_q=
A/A(h-q)A$, where $q$ is a complex scalar such that $1+q$ and $1+2q$
generate a free abelian subgroup of rank $2$ in $\CCx$. For such $q$, the
primitive ideals of $B_q$, as stated in \cite{\Vnc, Example 3.14},
consist of 
\roster
\item"$\bullet$" $\langle x_1-\alpha,x_2,x_3\rangle$,
$\langle x_1,x_2-\beta,x_3\rangle$, $\langle x_1,x_2,x_3-\gamma\rangle$,
for $\alpha,\beta,\gamma\in \CC$;
\item"$\bullet$" $\langle x_1\rangle$, $\langle
x_2\rangle$, $\langle x_3\rangle$, $\langle 0\rangle$.
\endroster
These can be computed by finding the prime ideals using \S9.6(a) and then
applying the Dixmier-Moeglin equivalence, which holds because $B_q$ is
a quantum affine space \cite{\qaff, Corollary 2.5}.

Observe that $\prim B_q$ is not homeomorphic to $\Pprim R$. For instance,
$\prim B_q$ has a generic point, while $\Pprim R$ does not. 

{\bf (c)} In contrast to the above, any generic $B_q$ is a member of a
flat family with a semiclassical limit $R'$ (the algebra $R$, but with a
different Poisson structure) such that $\prim B_q \approx
\Pprim R'$ and $\spec B_q \approx \Pspec R'$, by \cite{\GLet, Theorem
3.6}. 
\enddefinition

It would be very interesting to obtain criteria to
determine which flat families yield ``good'' semiclassical limits
relative to Conjecture 9.1. For quantum affine spaces and their Poisson
analogs, one good condition appears in the work of Oh, Park, and Shin
\cite{\OPS, Theorem 3.5} -- roughly, if the scalars appearing in the
defining Poisson brackets of the semiclassical limit arise from an
embedding into $k^+$ of the subgroup of
$\kx$ generated by the scalars appearing in the defining commutation
relations of the quantum affine space, then the prime and primitive
spectra of the quantum affine space are homeomorphic to the Poisson prime
and Poisson-primitive spectra of the semiclassical limit.

We close with an example of the ``simplest possible'' quantum group for
which primitive ideals match symplectic cores but not symplectic leaves.
There is no nontrivial multiparameter version of quantum $SL_2$, and to
deal with
$\O_{q,p}(SL_3(\CC))$ would require investigating $36$ families of
primitive ideals (indexed by $S_3\times S_3$, as in \cite{\HLT,
Corollary 4.5}). Instead, we look at a multiparameter quantization of
$GL_2$. It is convenient to use Takeuchi's original
presentation \cite{\Taktwo}.

\definition{9.10\. Example} For the classification of primitive ideals, we
assume only that $k$ is algebraically closed, and we choose a generic
pair of parameters $p,q\in\kx$, meaning that they generate a free abelian
subgroup of rank $2$ in $\kx$. We  restrict to $k=\CC$ and special
choices of $p$ and $q$ when setting up a semiclassical limit and
discussing symplectic leaves.

{\bf (a)} Define the two-parameter quantum $2\times2$ matrix algebra
$M_{q^{-1},p}$ as in \cite{\Taktwo}. This is the $k$-algebra with
generators $X_{11}$, $X_{12}$,
$X_{21}$, $X_{22}$ and relations
$$\xalignat2
X_{11}X_{12} &= qX_{12}X_{11}  &X_{11}X_{21} &=
p^{-1}X_{21}X_{11}\\
X_{21}X_{22} &= qX_{22}X_{21}  &X_{12}X_{22} &= p^{-1}X_{22}X_{12}  \\
X_{12}X_{21} &= (pq)^{-1}X_{21}X_{12}  &X_{11}X_{22}-X_{22}X_{11} &=
(q-p) X_{12}X_{21}\,.
\endxalignat$$
The element $D= X_{11}X_{22}- qX_{12}X_{21}$ is the quantum determinant
in $M_{q^{-1},p}$, but it is normal rather than central:
$$X_{ij}D = (pq)^{i-j}DX_{ij}$$
for all $i$, $j$ \cite{\Taktwo, \S2}. Since the powers of $D$ form an Ore
set, we can construct the Ore localization $A= A_{q^{-1},p}=
M_{q^{-1},p}[D^{-1}]$. There is a Hopf algebra structure on $A$
\cite{\Taktwo, \S2}, but we do not need that here.

For comparison with other presentations of multiparameter quantized
coordinate rings, we point out that
$A= \O_{pq^{-1},\bfp}(GL_2(k))$ in the  notation of \cite{\GMur, \S1.3;
\BrGd, \S I.2.4}), where $\bfp= \left[\smallmatrix 1&q^{-1}\\ q&1
\endsmallmatrix\right]$. In particular,
\cite{\BrGd, Corollary II.6.10} applies, implying that all prime ideals of
$A$ are completely prime.

Observe that $X_{12}$ and $X_{21}$ are normal in $A$, and so we can
localize with respect to their powers. Although $X_{11}$ is not normal,
its powers also form an Ore set (e.g., verify this first in
$M_{q^{-1},p}$, which is an iterated skew polynomial ring over
$k[X_{11}]$). Note that any ideal $I$ of $A$ which contains
$X_{11}$ also contains
$X_{12}X_{21}$, whence $D\in I$ and $I=A$. Hence, no prime ideal of $A$
contains
$X_{11}$, which means that no prime ideals of $A$ are lost in passing
from $A$ to the localization $A[X_{11}^{-1}]$.

{\bf (b)} The quotient $A/\langle X_{12},\,X_{21}\rangle$ is isomorphic
to a commutative Laurent polynomial ring
$k[x_{11}^{\pm1},x_{22}^{\pm1}]$. Hence, we know the prime ideals of $A$
containing $\langle X_{12},\,X_{21}\rangle$. The others correspond to
prime ideals in the localizations $\bigl( A/\langle X_{12}\rangle \bigr)
[X_{21}^{-1}]$, $\bigl( A/\langle X_{21}\rangle \bigr)
[X_{12}^{-1}]$, and $A[X_{11}^{-1},X_{12}^{-1},X_{21}^{-1}]$. We claim
that these localizations are simple algebras, from which it will follow
that the only prime ideals of $A$ not containing $\langle
X_{12},\,X_{21}\rangle$ are $\langle X_{12}\rangle$, $\langle
X_{21}\rangle$, and $\langle 0\rangle$.

First, $\bigl( A/\langle X_{12}\rangle \bigr)
[X_{21}^{-1}]$ is isomorphic to the algebra
$$T_1 := k\langle x^{\pm1},y^{\pm1},z^{\pm1} \mid xy=p^{-1}yx,\ xz=xz,\
yz= qzy \rangle.$$
Via \S9.6(a), we compute that $Z(T_1)=k$,
whence $T_1$ is simple.
Similarly, $\bigl( A/\langle X_{21}\rangle \bigr)
[X_{12}^{-1}]$ is simple.

Third, observe that $X_{22}= X_{11}^{-1}(D+qX_{12}X_{21})$ in
$A[X_{11}^{-1}]$, and so this algebra can be generated by
$X_{11}^{\pm1}$, $X_{12}$, $X_{21}$, $D^{\pm1}$. Consequently,
$A[X_{11}^{-1},X_{12}^{-1},X_{21}^{-1}]$ is isomorphic to the $k$-algebra
$T_3$ with generators $x^{\pm1}$, $y^{\pm1}$, $z^{\pm1}$, $w^{\pm1}$ and
relations
$$\xalignat3 xy &=qyx  &xz &= p^{-1}zx  &xw &= wx  \\
yz &= (pq)^{-1}zy  &yw &= (pq)^{-1}wy  &zw &= pqwz.
\endxalignat$$
Another application of \S9.6(a) shows that $T_3$ is simple,
establishing the claim.

Therefore, the prime ideals of $A$ consist of
\roster
\item"$\bullet$" the maximal ideals $\langle X_{11}-\lambda,\, X_{12},\,
X_{21},\, X_{22}-\mu \rangle$, for $\lambda,\mu\in \kx$;
\item"$(\diamondsuit)\ \bullet$" the ideals $\langle X_{12},\, X_{21},\,
f(X_{11},X_{22})
\rangle$, for irreducible polynomials $f(s,t) \in k[s^{\pm1},t^{\pm1}]$;
\item"$\bullet$" the ideals $\langle X_{12},\,X_{21}\rangle$, $\langle
X_{12}\rangle$, $\langle X_{21}\rangle$, and $\langle 0\rangle$.
\endroster

{\bf (c)} The torus $H= (\kx)^4$ acts on $A$ by $k$-algebra automorphisms
such that 
$$(\alpha_1,\alpha_2,\beta_1,\beta_2).X_{ij}= \alpha_i\beta_j X_{ij} 
\tag 9.10c$$ 
for all $i$, $j$. Only four of the prime ideals of $A$ are
$H$-stable, and thus \cite{\GLDM, Corollary 2.7(ii), Remark 5.9(i)}
implies that $A$ satisfies the Dixmier-Moeglin equivalence (cf\.
\cite{\BrGd, Corollary II.8.5(c)}). Therefore, the primitive ideals of
$A$ are
\roster
\item"$\bullet$" the maximal ideals $\langle X_{11}-\lambda,\, X_{12},\,
X_{21},\, X_{22}-\mu \rangle$, for $\lambda,\mu\in \kx$;
\item"$\bullet$" the ideals $\langle
X_{12}\rangle$, $\langle X_{21}\rangle$, and $\langle 0\rangle$.
\endroster

{\bf (d)} Now restrict to $k=\CC$, choose $\alpha\in \RR\setminus\QQ$,
assume that $q$ is transcendental over $\QQ(\alpha)$, and take $p=
1+\alpha(q-1)$. The assumptions on $\alpha$ and $q$ ensure that the
subgroup $\langle p,q\rangle \subseteq \CCx$ is free abelian of rank $2$,
as needed above. Our choice of $p$ is a first-order Taylor approximation
of $q^{\alpha}$, which is convenient for extension to polynomial rings.

Choose a Laurent polynomial ring $k[z^{\pm1}]$, set $z_\alpha=
1+\alpha(z-1)$, and let $B= M_{z^{-1},z_{\alpha}}$ over
$k[z^{\pm1},z_\alpha^{-1}]$ in the notation of \cite{\Taktwo, \S2}.
Thus, $B$ is the $k[z^{\pm1},z_\alpha^{-1}]$-algebra given by generators
$X_{11}$,
$X_{12}$,
$X_{21}$, $X_{22}$ and relations
$$\xalignat2
X_{11}X_{12} &= zX_{12}X_{11}  &X_{11}X_{21} &=
z_\alpha^{-1}X_{21}X_{11}\\
X_{21}X_{22} &= zX_{22}X_{21}  &X_{12}X_{22} &=
z_\alpha^{-1}X_{22}X_{12}  \\  
X_{12}X_{21} &= (zz_\alpha)^{-1}X_{21}X_{12}  &X_{11}X_{22}-X_{22}X_{11}
&= (z-z_\alpha) X_{12}X_{21}\,.
\endxalignat$$
Observe that $B$ is an iterated skew polynomial algebra over
$k[z^{\pm1},z_\alpha^{-1}]$, and so it is torsionfree over
$k[z^{\pm1}]$. This algebra has been arranged so that $B/(z-q)B \cong
M_{q^{-1},p}$ and $B/(z-1)B \cong \O(M_2(k))$. In $B$, the quantum
determinant is $D= X_{11}X_{22}- zX_{12}X_{21}$, and it is normal. We set
$C= B[D^{-1}]$ and observe that $C$ is a torsionfree $k[z^{\pm1}]$-algebra
such that $C/(z-q)C \cong A$ and $C/(z-1)C\cong R:= \O(GL_2(k))$.

Thus, $A$ is one of the quantizations of $R$ in the family of algebras
$C/(z-\gamma)C$. The semiclassical limit of this family is the algebra
$R$, equipped with the Poisson bracket determined by
$$\xalignat2
\{X_{11},X_{12}\} &= X_{11}X_{12}  &\{X_{11},X_{21}\} &= -\alpha
X_{11}X_{21}  \\
\{X_{21},X_{22}\} &= X_{21}X_{22}  &\{X_{12},X_{22}\} &= -\alpha
X_{12}X_{22}  \\
\{X_{12},X_{21}\} &= -(1+\alpha)X_{12}X_{21}  &\{X_{11},X_{22}\} &=
(1-\alpha)X_{12}X_{21} \,.
\endxalignat$$

To find the Poisson prime ideals of $R$, we can proceed in parallel with
part (b) above, using \S9.6(b) in place of \S9.6(a). We compute that the
Poisson prime ideals of
$R$ can be listed exactly as in
$(\diamondsuit)$. This yields an obvious bijection $\phi: \spec A
\rightarrow \Pspec R$ given by ``preservation of notation''. Clearly
$\phi$ and $\phi^{-1}$ preserve inclusions, and so $\phi$ is a
homeomorphism by Lemma 9.4(a). (Alternatively, one can easily identify the
closed sets in $\spec A$ and $\Pspec R$ and then check that $\phi$ and
$\phi^{-1}$ are closed maps.)

The torus $H$ acts on $R$ by Poisson algebra automorphisms satisfying
(9.10c), and only four Poisson prime ideals of $R$ are stable under
this action. Consequently, \cite{\GAth, Theorem 4.3} implies that $R$
satisfies the Poisson Dixmier-Moeglin equivalence. Thus, the Poisson-primitive ideals of $R$ are
\roster
\item"$(\spadesuit)\ \bullet$" the maximal ideals $\langle
X_{11}-\lambda,\, X_{12},\, X_{21},\, X_{22}-\mu \rangle$, for
$\lambda,\mu\in \kx$;
\item"$\bullet$" the ideals $\langle
X_{12}\rangle$, $\langle X_{21}\rangle$, and $\langle 0\rangle$,
\endroster
and therefore $\phi$ restricts to a homeomorphism $\prim A \rightarrow
\Pprim R$.

{\bf (e)} In view of $(\spadesuit)$, we can now
identify the symplectic cores in $GL_2(\CC) \approx \maxspec R$ with
respect to the Poisson structure under discussion. They are
\roster
\item "$\bullet$" the singletons $\left\{ \leftmat \lambda&0\\ 0&\mu
\rightmat \right\}$, for $\lambda,\mu \in\CCx$;  \smallskip
\item "$\bullet$" the sets $\leftmat \CCx&\CCx\\ 0&\CCx \rightmat$,
$\leftmat \CCx&0\\ \CCx&\CCx \rightmat$ and $\left\{ \leftmat
\lambda&\beta\\ \gamma&\mu
\rightmat \in GL_2(\CC) \biggm| \beta,\gamma\ne 0 \right\}$.
\endroster
The space of symplectic cores in $GL_2(\CC)$ is homeomorphic to $\Pprim
R$ by Lemma 9.3.

Since $\leftmat \CCx&\CCx\\ 0&\CCx \rightmat$ and
$\leftmat \CCx&0\\ \CCx&\CCx \rightmat$ are complex manifolds of odd
dimension, they cannot be symplectic leaves. In fact, each is the union
of a one-parameter family of symplectic leaves, which can be calculated
as in \cite{\GLet, Example 3.10(v)}. For instance, the symplectic leaves
contained in $\leftmat \CCx&\CCx\\ 0&\CCx \rightmat$ are the surfaces
$$\left\{ \leftmat \lambda&\beta\\ 0&\delta\lambda^\alpha \rightmat
\biggm| \lambda,\beta \in\CCx \right\},$$
for $\delta\in\CCx$.
\enddefinition

\head Acknowledgements \endhead

We thank J. Alev, K. A. Brown, I. Gordon, S. Kolb, U. Kraehmer, S.
Launois, T. H. Lenagan, E. S. Letzter, M. Lorenz, M. Martino, L. Rigal, M.
Yakimov, and J. J. Zhang for many discussions on the topics of this paper.


\Refs \widestnumber\key{{\bf 99}}

\ref\no\AvMV \by M. Adler, P. van Moerbeke, and P. Vanhaecke \book
Algebraic Integrability, Painlev\'e Geometry and Lie Algebras \publaddr
Berlin \yr 2004 \publ Springer-Verlag \endref

\ref\no\BGR \by W. Borho, P. Gabriel, and R. Rentschler \book Primideale
in Einh\"ullenden aufl\"osbarer Lie-Algebren \bookinfo Lecture Notes in
Math. 357 \publaddr Berlin \yr 1973 \publ Springer-Verlag \endref

\ref\no\Bro \by K. A. Brown \finalinfo Personal communication, Nov. 2008 \endref

\ref\no\BrGd \by K. A. Brown and K. R. Goodearl \book Lectures on
Algebraic Quantum Groups
\bookinfo Advanced Courses in Math. CRM Barcelona \publ Birkh\"auser
\publaddr Basel \yr 2002
\endref

\ref\no\BGY \by K. A. Brown, K. R. Goodearl, and M. Yakimov \paper
Poisson structures on affine spaces and flag varieties. I. Matrix affine
Poisson space \jour Advances in Math. \vol 206 \yr 2006 \pages 567-629
\endref

\ref\no\BrGr \by K.A. Brown and I. Gordon \paper Poisson orders,
representation theory, and symplectic reflection algebras \jour J. reine
angew. Math. \vol 559 \yr 2003 \pages 193-216 \endref

\ref\no\CG \by N. Chriss and V. Ginzburg \book Representation Theory and
Complex Geometry \publaddr Boston \yr 1997 \publ Birkh\"auser \endref

\ref\no\CoMc \by D. H. Collingwood and W. M. McGovern \book Nilpotent
orbits in semisimple Lie algebras \publ Van Nostrand Reinhold \publaddr
New York \yr 1993 \endref

\ref\no\Dix \by J. Dixmier \paper Repr\'esentations irr\'eductibles des
alg\`ebres de Lie r\'esolubles \jour J. Math. pures appl. \vol 45 \yr
1966 \pages 1-66 \endref

\ref\no \Dbook \bysame \book Enveloping Algebras \publ North-Holland
\publaddr Amsterdam
\yr 1977 \endref

\ref\no\FL \by D. R. Farkas and G. Letzter \paper Ring theory from
symplectic geometry \jour J. Pure Appl. Algebra \vol 125 \yr 1998 \pages
155-190 \endref

\ref\no\GMur \by K. R. Goodearl \paper  \paper Prime spectra of
quantized coordinate rings \inbook in Interactions between Ring Theory
and Representations of Algebras (Murcia 1998) \eds F.
Van Oystaeyen and M. Saor\'\i n \publaddr New York \yr 2000
\publ Dekker \pages 205-237  \endref

\ref\no\GAth \bysame \paper A Dixmier-Moeglin equivalence for
Poisson algebras with torus actions \paperinfo in Algebra and Its
Applications (Athens, Ohio, 2005) (D. V. Huynh, S. K. Jain, and S. R.
L\'opez-Permouth, Eds.) \jour Contemp. Math. \vol 419 \yr 2006 \pages
131-154 \endref

\ref\no\GLau \by K. R. Goodearl and S. Launois \paper The Dixmier-Moeglin
equivalence and a Gel'fand-Kirillov problem for Poisson polynomial
algebras \finalinfo preprint 2007, posted at
arxiv.org/abs/0705.3486 \endref

\ref\no\qaff \by K. R. Goodearl and E. S. Letzter \paper Prime and
primitive spectra of multiparameter quantum affine spaces \paperinfo in
Trends in Ring Theory. Proc. Miskolc Conf. 1996 (V. Dlab and L. M\'arki,
eds.) \jour Canad. Math. Soc. Conf. Proc. Series \vol 22
\yr 1998 \pages 39-58 \endref

\ref\no\GLDM \bysame \paper The
Dixmier-Moeglin equivalence in quantum coordinate rings and quantized
Weyl algebras \jour Trans. Amer. Math. Soc. \vol 352 \yr 2000 \pages
1381-1403 \endref

\ref\no\GLquo \bysame \paper Quantum $n$-space as a quotient of classical
$n$-space
\jour Trans. Amer. Math. Soc. \vol 352 \yr 2000 \pages 5855-5876  \endref

\ref\no\GLet \bysame \paper Semiclassical limits of quantum affine spaces
\jour Proc. Edinburgh Math. Soc. \finalinfo (to appear); posted at
arxiv.org/abs/0708.1091
\endref

\ref\no\GoYa \by K. R. Goodearl and M. Yakimov \paper Poisson structures
on affine spaces and flag varieties. II \jour Trans. Amer. Math. Soc.
\finalinfo (in press), posted at arxiv.org/abs/math.QA/0509075
\endref

\ref\no\Hay \by T. Hayashi \paper Quantum deformations of classical
groups\jour Publ. RIMS\vol 28\yr 1992\pages 57-81 \endref

\ref \no\Hod \by T. J. Hodges \paper Quantum tori and Poisson tori
\paperinfo Unpublished Notes, 1994 \endref

\ref \no\HoLeA \by T. J. Hodges and T. Levasseur\paper Primitive ideals of
$\bold C_q[SL(3)]$\jour Comm. Math. Phys.\vol 156\yr 1993\pages
581-605 \endref

\ref \no\HoLeB \bysame \paper Primitive ideals of
$\bold C_q[SL(n)]$\jour J. Algebra\vol 168\yr 1994\pages 455-468 \endref

\ref\no\HLT \by T. J. Hodges, T. Levasseur, and M. Toro\paper Algebraic
structure of multi-parameter quantum groups\jour Advances in Math.
\vol 126
\yr 1997
\pages 52-92 \endref

\ref\no\Hrt \by K. L. Horton \paper The prime and primitive spectra of
multiparameter quantum symplectic and Euclidean spaces \jour Communic. in
Algebra \vol 31 \yr 2003 \pages 2713-2743 \endref

\ref\no\Ing \by C. Ingalls \paper Quantum toric varieties \finalinfo
preprint, posted at kappa.math.unb.ca/\%7Ecolin/
\endref

\ref\no\Jos \by A. Joseph \paper On the prime and primitive spectra of
the algebra of functions on a quantum group\jour J. Algebra\vol 169\yr
1994\pages 441-511 \endref

\ref\no\Kam \by A. Kamita \paper Quantum deformations of 
certain prehomogeneous vector spaces III \jour Hiroshima Math. J. 
\vol 30 \yr 2000 \pages 79-115 \endref

\ref\no\Kirnilp \by A. A. Kirillov \paper Unitary representations of
nilpotent Lie groups \jour Russian Math. Surveys \vol 17 \yr 1962 \pages
57-110 \endref

\ref\no\Kirsurv \bysame \paper Merits and demerits of the
orbit method \jour Bull. Amer. Math. Soc. \vol 36 \yr 1999 \pages
433-488  \endref

\ref\no\KirC \bysame \book Lectures on the Orbit Method \bookinfo Grad.
Studies in Math. 64 \publaddr Providence \yr 2004 \publ American Math.
Soc. \endref

\ref\no\KlSc \by A. U. Klimyk and K. Schm\"udgen \book Quantum Groups and
their Representations \publaddr Berlin \yr 1997 \publ Springer-Verlag
\endref

\ref\no\KoSo \by L. I. Korogodski and Ya. S. Soibelman \book Algebras of
Functions on Quantum Groups: Part I \bookinfo Math. Surveys and
Monographs 56 \publ Amer. Math. Soc. \publaddr Providence \yr 1998
\endref

\ref\no\Mar \by M. Martino \paper The associated variety of a Poisson
prime ideal \jour J. London Math. Soc. (2) \vol 72 \yr 2005 \pages
110-120 \endref

\ref\no\Mat \by O. Mathieu \paper Bicontinuity of the Dixmier map \jour
J. Amer. Math. Soc. \vol 4 \yr 1991 \pages 837-863 \endref

\ref\no\McPe \by J. C. McConnell and J. J. Pettit \paper Crossed products
and multiplicative analogs of Weyl algebras \jour J. London Math. Soc.
(2) \vol 38 \yr 1988 \pages 47-55 \endref

\ref\no\Nou  \by M. Noumi\paper Macdonald's symmetric polynomials as
zonal spherical functions on some quantum homogeneous spaces 
\jour Advances in Math. \vol 123 \yr 1996 \pages 16-77 \endref

\ref\no\Ohone  \by S.-Q. Oh 
\paper Catenarity in a class of iterated skew polynomial rings
\jour Comm. Algebra \yr 1997 \vol 25(1) \pages 37-49
\endref

\ref\no\Ohtwo \bysame \paper Symplectic ideals of Poisson algebras
and the Poisson structure associated to quantum matrices \jour Comm.
Algebra
\vol 27 \yr 1999 \pages 2163-2180 \endref

\ref\no\Ohthree \bysame \paper Quantum and Poisson structures of multi-parameter symplectic and Euclidean spaces \jour J. Algebra \vol  319 \yr 2008 \pages 4485-4535 \endref

\ref\no\OPS \by S.-Q. Oh, C.-G. Park, and Y.-Y. Shin \paper Quantum
  $n$-space and Poisson $n$-space \jour Comm. Algebra \vol 30 \yr
2002\pages  4197-4209 \endref

\ref\no\Pol \by A. Polishchuk \paper Algebraic geometry of Poisson
brackets
\jour J. Math. Sci. \vol 94 \yr 1997
\pages 1413-1444 \endref

\ref\no\RTF \by N. Yu. Reshetikhin, L. A. Takhtadzhyan, and L. D.
Faddeev\paper Quantization of Lie groups and Lie algebras\jour Leningrad
Math. J.\vol 1\yr 1990\pages 193-225 \endref

\ref\no\Shaf \by I. Shafarevich \book Basic Algebraic Geometry \publaddr
Berlin \yr 1974 \publ Springer-Verlag \endref

\ref\no\SoiA \by Ya. S. Soibel'man \paper Irreducible representations of
the function algebra on the quantum group $SU(n)$, and Schubert
cells\jour Soviet Math. Dokl.\vol 40\yr 1990\pages 34-38  \endref

\ref\no\SoiB \bysame \paper The algebra of functions on a compact quantum
group, and its representations\jour Leningrad Math. J.\vol 2\yr
1991\pages 161-178 \moreref \paper Correction \paperinfo (Russian)\jour
Algebra i Analiz \vol 2
\yr 1990
\page 256 \endref

\ref\no\Str  \by E. Strickland \paper Classical invariant theory for 
the quantum symplectic group \jour Advances in Math. \vol 123 \yr 1996
\pages 78-90
\endref

\ref\no\Tak \by  M. Takeuchi \paper Quantum orthogonal and symplectic
groups and their embedding into quantum $GL$\jour Proc. Japan Acad.\vol
65 \yr 1989\pages 55-58 \endref

\ref \no\Taktwo \bysame \paper A two-parameter quantization of
$GL(n)$ {\rm(}summary{\rm)}\jour Proc. Japan. Acad.\vol 66 (A)\yr
1990\pages 112-114
\endref

\ref\no\TY \by P. Tauvel and R.W.T. Yu \book Lie Algebras and Algebraic
Groups \publaddr Berlin \yr 2005 \publ Springer-Verlag \endref

\ref\no\VS \by L. L. Vaksman and Ya. S. Soibel'man \paper Algebra of
functions on the quantum group $SU(2)$ \jour Func. Anal. Applic. \vol 22
\yr 1988 
pages 170-181 \endref

\ref\no\Vnc \by M. Vancliff \paper Primitive and Poisson spectra of
twists of polynomial rings \jour Algebras and Representation Theory \vol
2 \yr 1999 \pages 269-285 \endref

\ref\no\Van \by P. Vanhaecke \book Integrable Systems in the Realm of
Algebraic Geometry \bookinfo Lecture Notes in Math. 1638 \publaddr Berlin
\yr 1996 \publ Springer-Verlag \endref

\ref\no\Wein \by A. Weinstein \paper The local structure of Poisson
manifolds
\jour J. Diff. Geom. \vol 18 \yr 1983 \pages 523-557 \endref

\endRefs
\enddocument